\documentclass[10 pt, journal, final, twocolumn]{IEEEtran}

\usepackage{enumerate}
\usepackage{amsmath, amsthm, amssymb}
\usepackage{float}
\usepackage{subcaption}
\usepackage{flexisym}
\usepackage{mathtools}
\usepackage{breqn}
\usepackage{array}
\usepackage{tikz}
\usepackage{cite}
\usetikzlibrary{shapes.geometric, arrows}
\usepackage[font = small]{caption}
\usepackage{paralist}
\usepackage{mathrsfs}
\usepackage{algorithm}
\usepackage{algpseudocode}
\allowdisplaybreaks[1]

\newtheorem{lemma}{Lemma}
\newtheorem{theorem}{Theorem}

\DeclareMathOperator*{\trace}{trace}
\DeclareMathOperator*{\diag}{diag}

\begin{document}
\title{Optimal Attack Strategies Subject to Detection Constraints Against Cyber-Physical Systems}

\author{Yuan Chen, Soummya Kar, and Jos\'{e} M. F. Moura
\thanks{Yuan Chen {\{(412)-268-7103\}}, Soummya Kar {\{(412)-268-8962\}}, and Jos\'{e} M.F. Moura {\{(412)-268-6341, fax: (412)-268-3890\}} are with the Department of Electrical and Computer Engineering, Carnegie Mellon University, Pittsburgh, PA 15217 {\tt\small \{yuanche1, soummyak, moura\}@andrew.cmu.edu}}
\thanks{This material is based upon work supported by the Department of Energy under Award Number DE-OE0000779 and by DARPA under agreement number DARPA FA8750-12-2-0291. The U.S. Government is authorized to reproduce and distribute reprints for Governmental purposes notwithstanding any copyright notation thereon. The views and conclusions contained herein are those of the authors and should not be interpreted as necessarily representing the official policies or endorsements, either expressed or implied, of DARPA or the U.S. Government.}}
\maketitle

\begin{abstract}
	This paper studies an attacker against a cyber-physical system (CPS) whose goal is to move the state of a CPS to a target state while ensuring that his or her probability of being detected does not exceed a given bound. The attacker's probability of being detected is related to the nonnegative bias induced by his or her attack on the CPS's detection statistic. We formulate a linear quadratic cost function that captures the attacker's control goal and establish constraints on the induced bias that reflect the attacker's detection-avoidance objectives. When the attacker is constrained to be detected at the false-alarm rate of the detector, we show that the optimal attack strategy reduces to a linear feedback of the attacker's state estimate. In the case that the attacker's bias is upper bounded by a positive constant, we provide two algorithms -- an optimal algorithm and a sub-optimal, less computationally intensive algorithm -- to find suitable attack sequences. Finally, we illustrate our attack strategies in numerical examples based on a remotely-controlled helicopter under attack. 
\end{abstract}

\section{Introduction}\label{sect: introduction}
Security vulnerabilities in cyber-physical systems (CPS), systems that interface sensing, communication, and control with an underlying physical process, allow for sophisticated cyber attacks that cause catastrophic physical harm. In the past, events such as StuxNet~\cite{Cardenas} and the Maroochy Sewage Control Incident~\cite{CardenasOld} have demonstrated the vulnerability of industrial processes. More recently, cyber-physical attacks have targeted automobiles~\cite{CarAttack}, military vehicles~\cite{Iran}, and commercial drones~\cite{DroneHack}. These examples show that CPS remain susceptible to cyber-attacks, and, in response, there have been significant efforts to improve the security of CPS.

Part of the effort in improving cyber-physical security has been devoted to categorizing different types of attacks and developing security countermeasures for each type~\cite{Cardenas, TeixeiraModels}. One particular type of attack is the integrity attack, in which an attacker manipulates the CPS's sensor readings and alters its actuator control signals~\cite{TeixeiraModels, Mo}. Prior work has analyzed CPS security sensor attacks, determining the fundamental limits of attack detection~\cite{Liu} and developing methods to reconstruct sensor attacks~\cite{Fawzi, Shoukry, Pajic}. Existing work has also studied the capabilities of an integrity attacker, relating the ability of the attacker to perform undetectable attacks to certain geometric control-theoretic properties of the CPS~\cite{Pasqualetti, ChenICASSP, ChenSideInfo}. For systems affected by process and sensor noise, references~\cite{Mo} and~\cite{MoWorkshop} characterize the state estimation error caused by an attacker who tries to avoid detection.

In addition to analyzing the \textit{ability} of an integrity attacker to cause damage and evade detection, prior work has also studied \textit{how} an attacker should behave in order to achieve his or her objectives. Reference~\cite{ScalarAttack} considers a noisy CPS and designs an attack to optimally disrupt the system's feedback controller, while avoiding detection. Instead of attackers who seek to cause general disruption and damage to a CPS, our previous work studies attackers with specific control objectives~\cite{ChenACC, ChenAttack1}. In~\cite{ChenAttack1}, we considered an attacker whose goal is move the CPS to a target state while evading detection, formulated a cost function that penalized the deviation from the target state and the magnitude of the detection statistic, and we determined the optimal attack for such a cost function reduced to a linear feedback of the attacker's state estimate.

This paper studies an attacker who wishes to move the system to a target state, but, unlike~\cite{ChenAttack1}, we impose an explicit bound on the probability of an attack being detected. We model the CPS as a linear dynamical system subject to sensor and process noise equipped with a Kalman filter for state estimation and an LQG controller. The CPS uses a $\chi^2$ detector as an attack detector, which reports an attack if the energy of the Kalman filter innovation exceeds a certain threshold~\cite{Willsky}.
{\color{black}This model has been used in the literature to model CPS under attack (see, e.g.,~\cite{MoAuthentication, MoScada}).} The attacker's goal is to design a sequence of attacks that counters the system's LQG controller and minimizes the deviation of the CPS's state from the target state subject to a bound on the (non-negative) bias induced on the $\chi^2$ detection statistic.

We define a linear quadratic cost function that captures the attacker's control objective by penalizing the distance between the system's state and the target state. Then, we formulate the attack design problem as an optimization problem of finding a sequence of attacks that minimizes the cost function subject to an upper bound constraint on the bias induced in the detection statistic. {\color{black}This differs from our previous work~\cite{ChenAttack1} that studies \textit{unconstrained} attack design where the attacker's detection avoidance goals is an additional term of the overall cost function. This paper, unlike~\cite{ChenAttack1}, requires attacks at each time step to satisfy explicit constraints. Since we compute the optimal attack sequence in a causal manner, we must ensure that attacks at each time step are  \textit{recursively feasible}~\cite{RecursiveFeasibility} to guarantee that it is possible to satisfy the constraints of future time steps. We use geometric control properties (similar to those studied in~\cite{Rappaport}) of the CPS model to express, the constraint placed on the detection statistic bias as a linear constraint on the attack at each time step. From a practical perspective, this paper provides guarantees on the optimal attacks' probability of being detected (reference~\cite{ChenAttack1} does not provide such guarantees).} 

{\color{black} We consider separately two cases:
\begin{inparaenum}[1)] 
	\item when the induced bias is constrained to be zero and 
	\item in which the induced bias is upper bounded by a positive constant. 
\end{inparaenum}
When the bias is zero, which restricts the attacker to be detected at the false alarm rate of the detector, we apply constrained dynamic programming to show that the optimal attack reduces to a linear feedback strategy. For bounded bias, we provide two algorithms to determine a suitable sequence of attacks. The first algorithm is more computationally intensive but finds an optimal sequence of attacks. The second, less computationally intensive, algorithm finds a (sub-optimal) sequence of attacks that satisfies the detection constraint.  A preliminary version of part of this work appears in~\cite{ChenCDC}, designing attacks constrained to be detected at the false alarm rate, but when the CPS are \textit{not} equipped with LQG controllers.} {\color{black} When the CPS is equipped with its own controller, which we consider here, the attacker must account for system input in designing his or her attack.}

The rest of the paper is organized as follows. Section~\ref{sect: background} provides the model and assumptions for the CPS and attacker, reviews the $\chi^2$ detector and the concept of recursive feasibility, and formally states the problems we address. In Section~\ref{sect: feasibility}, we determine the set of all recursively feasible attacks at each time step. In Section~\ref{sect: equalityConstraints}, we use dynamic programming to find an optimal strategy when the attacker's probability of being detected is constrained to be the detector's false alarm rate. Section~\ref{sect: inequalityConstraints} studies the case when the bias induced in the detection statistic is upper bounded by a positive constant; we provide two algorithms for computing attack sequences that achieve the attacker's objectives. We provide numerical examples of a remotely-controlled helicopter under attack (from each of our proposed strategies) in Section~\ref{sect: examples}, and we conclude in Section~\ref{sect: conclusion}.

\section{Background}\label{sect: background}
\subsection{Notation}\label{sect: notation}
Let $\mathbb{R}$ denote the reals, $\mathbb{R}^n$ denote the space of $n$-dimensional real (column) vectors, and $\mathbb{R}^{m \times n}$ denote the space of real $m$ by $n$ matrices. The multivariate Gaussian distribution with mean $\mu$ and covariance $\Sigma$ is denoted as $\mathcal{N} \left( \mu, \Sigma \right)$. The $n$ by $n$ identity matrix is denoted as $I_n$. For a matrix $M$, $\mathscr{R} (M)$ denotes the range space of $M$, $\mathscr{N} (M)$ denotes the null space of $M$, and $M^\dagger$ denotes the Moore-Penrose pseudoinverse. For a symmetric matrix $S = S^T$, $S \succeq 0$ denotes that $S$ is positive semidefinite, and $S \succ 0$ denotes that $S$ is positive definite. For $S \succ 0$, let $\left \lVert \cdot \right \rVert_S$ denote the $S$-weighted $2$-norm. That is, for $S \in \mathbb{R}^{n\times n}$ and $x \in \mathbb{R}^{n}$,
$\left\lVert x\right\rVert_S = \left\lVert S^{\frac{1}{2}} x \right\rVert_2 = \sqrt{x^T S x}.$

\subsection{System Model}\label{sect: model}
We use the same, linearized\footnote{If the CPS is nonlinear, the model~\eqref{eqn: ssModel1} represents its dynamics after linearization about an operating point.} CPS model as~\cite{ChenAttack1}:
\begin{equation}\label{eqn: ssModel1}
\begin{split}
	x_{t+1} &= Ax_{t} + Bu_{t} + \Gamma e_t + w_t, \\
	y_t& = Cx_{t} + \Psi e_t + v_t,
\end{split}
\end{equation}
where $x_t \in \mathbb{R}^n$ describes the system's state, $u_t \in \mathbb{R}^m$ is the system input, $e_t \in \mathbb{R}^s$ is the attacker's input{\color{black}\footnote{\color{black}The attack $e_t$ also models the case in which the attacker may separately attack the CPS's actuators and sensors. Define $e_t = \left[\begin{array}{cc} {e_t^d}^T & {e_t^d}^T \end{array} \right]^T$, where $e_t^d$ and $e_t^s$ are the actuator and sensor attacks, respectively, and define $\Gamma = \left[\begin{array}{cc} \widetilde{\Gamma} & 0 \end{array} \right]$ and $\Psi = \left[\begin{array}{cc} 0 & \widetilde{\Psi} \end{array} \right].$}}, and $w_t$ and $v_t$ are the process and sensor noise, respectively. The sensor and process noise are independently, identically, distributed (i.i.d.) in time and mutually independent; $w_t$ has distribution $\mathcal{N}\left(0, \Sigma_w\right)$, and $v_t$ has distribution $\mathcal{N}\left( 0, \Sigma_v \right)$, with $\Sigma_w, \Sigma_v \succ 0$. The system starts running at time $t = -\infty$, and the initial state of the system $x_{-\infty}$ has distribution $\mathcal{N} \left( 0, \Sigma_x \right)$ with $\Sigma_x \succ 0$ and is independent of the noise processes. The pair $(A, C)$ is observable, and the pair $(A, B)$ is controllable. The matrices $\Gamma$ and $\Psi$ describe the attacker. The model~\eqref{eqn: ssModel1} is commonly adopted in studies of CPS under attack~\cite{ChenAttack1, MoScada, MoAuthentication}.

The system knows the matrices $A, B$, and $C$ and the statistics of the noise processes and initial state, but does not know the matrices $\Gamma$ and $\Psi$ (since they describe the attacker). The system causally knows the system input $u_t$ and the sensor output $y_t$, but not the attack $e_t$. We assume the system's goal is to regulate the system state to the origin. Because the system cannot directly observe the state $x_t$, it uses its sensor measurements $y_t$ to construct an estimate of the state using a Kalman filter. Then, the system performs feedback control on the state estimate to regulate the state to the origin. The system constructs its Kalman filter and controller assuming nominal operating conditions (i.e., $e_t = 0$ for all $t$). 

Under nominal operating conditions, the system's Kalman filter calculates $\widehat{x}_{t}$, the minimum mean square error (MMSE) estimate of $x_t$ given all sensor measurements up to time $t$ and input up to time $t-1$. Since the system starts at $t = -\infty$, the Kalman filter has fixed gain:
\begin{align}\label{eqn: LimitingKalman}
K &= PC^T(CPC^T + \Sigma_v)^{-1},\\
\begin{split}
P &= APA^T + \Sigma_w \\
	&\quad - APC^T (CPC^T + \Sigma_v)^{-1}CPA^T, \label{eqn: covarianceEquation}
\end{split}\\
\widehat{x}_t & = \widehat{x}_{t|t-1} + K\left(y_{t} - C \widehat{x}_{t|t-1}\right), \label{eqn: filterEquation}\\
\widehat{x}_{t+1|t} &=A\widehat{x}_t + Bu_t. \label{eqn: kalmanPrediction}
\end{align}
To regulate the state $x_t$, the system has a feedback controller of the form
\begin{equation}\label{eqn: sysController}
	u_t = L \widehat{x}_t,
\end{equation}
where the feedback matrix $L$ is chosen such that $A+BL$ is stable. One controller that takes the form of equation~\eqref{eqn: sysController} is the infinite horizon LQG controller that minimizes the cost function
	$J_{\text{CPS}} = \lim_{T \rightarrow \infty} \frac{1}{2T+1} \mathbb{E} \left[\sum_{t = -T}^{T} x_t^T Q^' x_t + u_t^T R^' u_t\right],$
where $Q^' \succeq 0$, $R^' \succ 0$, and the pair $(A, Q^')$ is observable. 

The CPS is equipped with a $\chi^2$ detector to determine if, for some $t$, $e_t \neq 0$. The $\chi^2$ attack detector~\cite{Willsky} uses the innovations sequence of the Kalman filter, $\nu_t$, defined as
	$\nu_t = y_t - C\widehat{x}_{t|t-1},$
to determine whether or not an attack has occurred. The term $\widehat{x}_{t|t-1}$ is the MMSE estimate of $x_t$ given all sensor measurements and system input up to time $t-1$, assuming nominal operating conditions. When there is no attack (i.e., $e_t = 0$ for all $t$), the innovations sequence is i.i.d. $\mathcal{N} \left( 0, \Sigma_\nu \right)$, where $\Sigma_\nu = CPC^T + {\color{black}\Sigma_v}$, and $\nu_t$ is orthogonal to $\widehat{x}_{t|t-1}$~\cite{Speyer}. {\color{black}The $\chi^2$ reports an attack if the statistic
	$g_t = \sum_{k = t-t^'+1}^t \nu_k^T \Sigma_\nu^{-1} \nu_k,$
where $S$ is the window size of the detector, exceeds a threshold $\tau$, which is chosen \`{a} priori to balance the false alarm and missed detection probabilities~\cite{Willsky}. In this paper, we consider a $\chi^2$ detector with window size $t^' =1$, so $g_t = \nu_t^T \Sigma_\nu^{-1}\nu_t$.} 

{\color{black} There are attack detectors other than the $\chi^2$ attack detector (see, e.g.,~\cite{Fawzi, Pasqualetti, PajicCDC}.) These detectors require noiselessness~\cite{Fawzi, Pasqualetti}, bounded energy noise~\cite{PajicCDC}, or batch measurements~\cite{Fawzi, PajicCDC}. This paper studies attacks against CPS under broader conditions on the noise and provides a recursive solution. We consider an on-line $\chi^2$ detector for systems with process and sensor noise.} The linear, state space model with a Kalman filter, feedback controller, and $\chi^2$ attack detector {\color{black} with window size $1$} is a standard model for a CPS subject to attack~\cite{ChenAttack1, MoScada, MoAuthentication}.

\subsection{Attacker Model}\label{sect: attacker}
The attacker knows the system model and statistical properties, the controller feedback matrix $L$, and which sensors and actuators he or she can attack (i.e. the matrices $\Gamma$ and $\Psi$){\color{black}\footnote{\color{black}In future work, we will consider defense strategies against such attackers and study the interaction between the attacker and CPS in a game theoretic framework. Thus, in this paper, we assume the worst-case, most powerful attacker who knows the system model perfectly. In addition, future work will study attack strategies that are robust to imperfections in the attacker's knowledge of the system model.}}. Following~\cite{ChenAttack1}, we assume, without loss of generality, that the matrix $\left[\begin{array}{c} \Gamma \\ \Psi \end{array} \right]$ is injective. The attacker causally knows the sensor output $y_t$ and the attack $e_t$. Additionally, the attacker causally knows {\color{black} $\widetilde{y}_t = Cx_t + v_t$}, the value of the sensor output at time $t$ before it is altered by the attack at time $t$.

The attacker performs Kalman filtering, separately from the system, to estimate the state. The attacker also uses his or her knowledge to compute the estimate produced by the \textit{system's} Kalman filter. The attacker knows $L$ and the system's state estimate, so he or she knows $u_t$. {\color{black} We design attack strategies that depend on the attacker knowing the system's input. In general, so long as the attacker knows $u_t$ for all $t$, the CPS's control input need not be restricted to the form of equation~\eqref{eqn: sysController}. For this paper, we only consider the case of feedback control, but our methodology may be tailored toward other control laws.} The attack begins at time $t=0$, i.e., for $t = -\infty, \dots, -1$, $e_t = 0$. During the time interval $t = -\infty, \dots, 0$, the attacker observes the system output and keeps track of the state estimate $\widehat{x}_t$. 

The attacker's objective is to design an attack sequence over the finite time interval $t=0$ to $t=N$ that moves the system state to a target state $x^*$ while satisfying a detection-avoidance constraint. The attacker chooses the sequence \[\gamma(0, N) = \left\{e_0, \dots, e_N\right\},\] to accomplish his or her goals such that at time $t$, the attack $e_t$ only depends on the attacker's available information at time $t$, $\mathcal{I}_t$. Following~\cite{ChenAttack1}, $\mathcal{I}_t$ is the classical information pattern~\cite{Speyer}:
	$\mathcal{I}_0 = \left\{\widetilde{y}_0\right\}, \mathcal{I}_{t+1} = \left\{\mathcal{I}_{t}, \widetilde{y}_{t+1}, e_{t}\right\}.$
If a nonzero attack occurs, the attacker's Kalman filter then produces a different estimate than the system's Kalman filter and becomes:{\color{black}\begin{align}\label{eqn: attackerPredictor}
	\widetilde{x}_{t+1|t} &= A\widehat{x}_t + Bu_t + \Gamma e_t, \\
	\widetilde{x}_t &= \widehat{x}_{t|t-1} + K\left(\widetilde{y}_t - C\widehat{x}_{t|t-1}\right), \label{eqn: attackerFilter}
\end{align}}The attacker's Kalman filter produces the MMSE state estimate given $\mathcal{I}_t$, i.e., $\widetilde{x}_t = \mathbb{E} \left[ x_t \vert \mathcal{I}_t \right]$. 

The attack $\gamma(0, t)$ induces a bias $\epsilon_t$ in the {\color{black} system's} innovation $\nu_t$. {\color{black} Under an attack $\gamma(0, t)$, we have
	$\nu_t = \nu_t^0 + \epsilon_t,$
where $\nu_t^0$ is the value of the system's innovation in the case that there had been no attack (i.e., $e_0 = e_1 = \dots = e_t = 0$).} 
The following state space dynamical system describes the relationship between $\gamma(0, t)$ and $\epsilon_t$~\cite{ChenAttack1}:
\begin{equation}\label{eqn: eEpsilonSystem}
\begin{split}
	\theta_{t+1}  &= \widehat{\mathcal{A}} \theta_t + \widehat{\mathcal{B}} e_t, \\
	\epsilon_t &= \widehat{\mathcal{C}} \theta_t + \widehat{\mathcal{D}} e_t,
\end{split}
\end{equation}
where $\widehat{\mathcal{A}} = \left[\begin{array}{ccc} (I_n-KC)A + BL & KC & KC \\ ABL & A & 0 \\ 0 & 0 & A \end{array}\right]$, $\widehat{\mathcal{B}}^T = \left[\begin{array}{ccc} \left(K\Psi\right)^T & 0^T & \Gamma^T \end{array}\right]$, $\widehat{\mathcal{C}} = \left[\begin{array}{ccc} -CA & C & C \end{array}\right]$, $\widehat{\mathcal{D}} = \Psi$, and $\theta_0 = 0$.

The $\Sigma_\nu^{-1}$-weighted $2$-norm of $\epsilon_t$ relates to the probability of the attack being detected at time $t$~\cite{MoDegradation}. Let $P_{D, t} = \mathbb{P} \left( g_t > \tau\right)$ be the detection probability at time $t$. {\color{black} If $\left\lVert \epsilon_t \right\rVert_{\Sigma_{\nu}^{-1}}^2 = 0$, then, for any positive detector window size, the probability of detection at time $t$ is equal to the false alarm probability of the $\chi^2$ detector, since there is no induced bias in $\nu_t$. For nonzero bias (and detector window size $1$), the following lemma relates the bound on $\left\lVert \epsilon_t \right\rVert_{\Sigma_{\nu}^{-1}}^2$ to the probability of being detected.}

\begin{lemma}[Detection Probability Bound~\cite{MoDegradation}]\label{thm: detectionBound}
	For any $\delta \in \left(0, \tau\right)$, if $\left\lVert \epsilon_t \right\rVert_{\Sigma_{\nu}^{-1}}^2 \leq \delta$, then
	\[ P_{D, t} \leq \mathbb{P}\left( g_t^0 > \left(\sqrt{\tau} - \sqrt{\delta} \right)^2 \right),\]
	where $g_t^0$ is the value of the statistic $g_t$ when there is no attack\footnote{The statistic $g_t^0$ is i.i.d. (in time) $\chi^2$ with $p$ degrees of freedom.}.
\end{lemma}

To model the attacker's control objectives, define the cost function:
\begin{equation}\label{eqn: costFunction}
	J = \mathbb{E} \left[ \sum_{t = 0}^N \left\lVert \left(x_t - x^*\right) \right \rVert_{Q_t}^2 \right],
\end{equation}
with $Q_t \succ 0$. The cost function $J$ penalizes deviation of the state from the target state. The attacker's goal is to design an attack $\gamma(0, N)$ that achieves cost
\begin{equation}\label{eqn: optimalCost}
	\begin{array}{ccl}	
	J^* = & \underset{\gamma(0, N)}{\text{min}} & \mathbb{E} \left[ \sum\limits_{t = 0}^N \left\lVert \left(x_t - x^*\right) \right \rVert_{Q_t}^2 \right] \\
		& \text{s.t.} & \left\lVert \epsilon_t \right \rVert_{\Sigma_\nu^{-1}}^2 \leq \delta, \forall t = 0, \dots, N
	\end{array},
\end{equation}
the minimum cost of $J$ subject to constraints on the $\Sigma_{\nu}^{-1}$-weighted $2$-norm of $\epsilon_t$.
The constraints in the optimization problem~\eqref{eqn: optimalCost} model the attacker's goals of evading detection.

\subsection{Recursive Feasibility}\label{sect: feasibilityBackground}
The attacker designs the attack in real time: at time $t$, the attacker chooses the attack $e_t$ based on his or her information $\mathcal{I}_t$. Note that the constraint $\left\lVert \epsilon_t \right \rVert_{\Sigma_\nu^{-1}}^2 \leq \delta$ in~\eqref{eqn: optimalCost} is for all times $t = 0, \dots, N$. It is necessary that the attack be recursively feasible~\cite{RecursiveFeasibility}: the attack $e_t$ must be chosen such that $\left\lVert \epsilon_t \right \rVert_{\Sigma_\nu^{-1}}^2 \leq \delta$ and, for all future times $t+1, \dots, N$, there exist attacks $e_{t+1}, \dots, e_N$ such that $\left\lVert \epsilon_{t+1} \right \rVert_{\Sigma_\nu^{-1}}^2 \leq \delta, \dots, \left\lVert \epsilon_N \right \rVert_{\Sigma_\nu^{-1}}^2 \leq \delta$. The recursive feasibility of~\eqref{eqn: optimalCost} is related to the output minimization problem presented in~\cite{Rappaport}:
\begin{lemma}[\!\!\cite{Rappaport}]\label{lem: outputRiccati}
	Consider the system in~\eqref{eqn: eEpsilonSystem} with arbitrary initial state $\theta_0$. Then, for any $k = 1, 2, \dots$,
	\begin{equation}\label{eqn: outputMinimization}
		\min_{e_0, \dots, e_{k-1}} \sum_{t = 0}^{k-1} \left\lVert \epsilon_t  \right\rVert_{\Sigma_{\nu}^{-1}}^2 = \theta_0^T \widehat{\mathbf{P}}_k \theta_0,
	\end{equation}
	where $\widehat{\mathbf{P}}_k$ follows the solution to the Riccati equation
	\begin{equation}\label{eqn: outputRiccati}
	\begin{split}
		&\widehat{\mathbf{P}}_{k+1} = \!\!\widehat{\mathcal{A}}^T \widehat{\mathbf{P}}_k \widehat{\mathcal{A}} + \widehat{\mathcal{C}}^T \Sigma_{\nu}^{-1} \widehat{\mathcal{C}} - \!\!\Big( \widehat{\mathcal{D}}^T \Sigma_{\nu}^{-1} \widehat{\mathcal{C}} +  \widehat{\mathcal{B}}^T \widehat{\mathbf{P}}_k \widehat{\mathcal{A}} \Big)^T\!\!\\
	&\times{\color{black}\left(\widehat{\mathcal{D}}^T \Sigma_{\nu}^{-1} \widehat{\mathcal{D}} + \widehat{\mathcal{B}}^T \widehat{\mathbf{P}}_k \widehat{\mathcal{B}} \right)^{\dagger}}\Big(\widehat{\mathcal{D}}^T \Sigma_{\nu}^{-1} \widehat{\mathcal{C}} + \widehat{\mathcal{B}}^T \widehat{\mathbf{P}}_k \widehat{\mathcal{A}} \Big),
	\end{split}
	\end{equation}
	with $\widehat{\mathbf{P}}_0 = 0$. 
\end{lemma}
\noindent Furthermore, the matrix $\widehat{\mathbf{P}}_k$ is positive semidefinite, and $\min_{e_0, \dots, e_{k-1}} \sum_{t = 0}^{k-1} \left\lVert \epsilon_t  \right\rVert_{\Sigma_{\nu}^{-1}}^2 = 0$ if and only if $\widehat{\mathbf{P}}_k \theta_0 = 0$~\cite{Rappaport}.

\subsection{Augmented State Space Notation}\label{sect: augmentedNotation}
For the remainder of this paper, we use the augmented state space description of the cyber-physical system and attacker provided in~\cite{ChenAttack1}. Define the augmented state 
\begin{equation}
	{\color{black}\xi_t^T = \left[ \begin{array}{cccc} \widetilde{x}_t^T & \theta_t^T & \left(\widehat{x}_t^0\right)^T & {x^*}^T \end{array} \right],}
\end{equation}
where $\widehat{x}_t^0$ denotes the system's state estimate in the case that $e_0 = \dots = e_t = 0$.  The state $\xi_t$ follows the dynamics
\begin{equation}\label{eqn: xiDynamics}
	\xi_{t+1} = \mathcal{A} \xi_t + \mathcal{B} e_t + \mathcal{K}{\color{black} \widetilde{\nu}_{t+1}},
\end{equation}
where {\color{black}$\widetilde\nu_{t+1}$} denotes the attacker's innovation at time $t+1$, ${\color{black}\mathcal{A} = \left[\begin{array}{cccc} A & BL\Omega & BL & 0 \\ 0 & \widehat{\mathcal{A}} & 0 & 0 \\ 0 & 0 & A+BL & 0 \\ 0 & 0 &0 &I_n\end{array} \right]}$, $\mathcal{B}^T = \left[\begin{array}{cccc} \widetilde{\mathcal{B}}^T & \widehat{\mathcal{B}}^T & 0^T & 0^T \end{array} \right], \mathcal{K}^T =  \left[\begin{array}{cccc} K^T & 0 & K^T &  0 \end{array}\right]^T$, $\Omega = \left[\begin{array}{ccc} (I_n-KC)A + BL & KC & KC \end{array}\right]$, and $\widetilde{\mathcal{B}} = \Gamma + BLK\Psi$.

Further define $\widetilde{\mathcal{C}} = \left[\begin{array}{cccc} 0 & \widehat{\mathcal{C}} & 0 & 0 \end{array}\right]$, $\widetilde{\mathcal{D}} = \Psi$, and $\mathcal{H} = \left[\begin{array}{cccc} I_n & 0 & 0 & -I_n \end{array}\right]$. 
Then, we have
\begin{align}
	\epsilon_t &= \widetilde{\mathcal{C}} \xi_t  + \widetilde{\mathcal{D}} e_t,\label{eqn: epsilonXi} \\
	\color{black} \widetilde{x}_t - x^* &= \mathcal{H}\xi_t.
\end{align}
One important property of $\xi_t$ is that, given $\mathcal{I}_t$, the attacker can exactly determine the value of $\xi_t$~\cite{ChenAttack1}. Accordingly, the attacker can use $\xi_t$ to determine his or her attack at time $t$. 

Following~\cite{ChenAttack1} and~\cite{Speyer}, we manipulate the cost function $J$ by substituting ${\color{black}x_t = \widetilde{x}_t + n_t}$, where $n_t$ is the estimation error. It is well known~\cite{Speyer} that, given $\mathcal{I}_t$, $n_t$ is conditionally {\color{black}distributed} as $\mathcal{N} \left(0, \widehat{P} \right)$, where $\widehat{P} = P - PC^T \Sigma^{-1}_{\nu}CP,$
and $n_t$ is conditionally orthogonal to $\widehat{x}_t$. Performing this substitution, the optimal attack design problem becomes.
\begin{equation}\label{eqn: optimalAttack}
	\begin{array}{cl}
		\underset{\gamma(0, N)} {\text{min}} &  \sum\limits_{t=0}^N \trace\left(\widehat{P} Q_t \right) + \mathbb{E} \left[ \sum\limits_{t=0}^N \left\lVert \mathcal{H} \xi_t \right \rVert_{Q_t}^2 \right] \\
		\text{s.t.} & \left\lVert \epsilon_t \right \rVert_{\Sigma_\nu^{-1}}^2 \leq \delta, \forall t = 0, \dots, N
	\end{array},
\end{equation}
where $\sum_{t=0}^N \trace\left(\widehat{P} Q_t \right)$ does not depend on $\gamma(0, N)$. 

\subsection{Problem Statement}\label{sect: problem}
This paper addresses three main problems. Consider the optimal attack design problem~\eqref{eqn: optimalAttack}.
First, determine, for any $\delta \geq 0$ and any time $t = 0, \dots, N$, the set of recursively feasible attacks. Second, find an optimal attack sequence $\gamma(0, N)$ when $\delta = 0$. This corresponds to finding the optimal attack under the constraint that the probability of being detected at any time $t$ is equal to the false alarm probability of the detector. 
Third, find an optimal attack sequence $\gamma(0, N)$ when $\delta > 0$. 

\section{Feasibility Sets}\label{sect: feasibility}
In this section, we determine which attacks $e_t$ are recursively feasible at time $t$. Recursively feasible attacks are attacks $e_t$ such that $\left\lVert \epsilon_t \right \rVert_{\Sigma_\nu^{-1}}^2 \leq \delta$ and there exists $e_{t+1}, \dots, e_{N}$ such that $\left\lVert \epsilon_{t+1} \right \rVert_{\Sigma_\nu^{-1}}^2 \leq \delta, \dots, \left\lVert \epsilon_N \right \rVert_{\Sigma_\nu^{-1}}^2 \leq \delta$. From equations~\eqref{eqn: xiDynamics} and~\eqref{eqn: epsilonXi}, we see that the recursively feasibility of an attack $e_t$ depends on the state $\xi_t$. Define the sets $\Xi_t$, $t = 0, \dots, N$ as follows:
\begin{equation}\label{eqn: XiDef}
	\begin{split}
		\Xi_N &= \left \{ \xi_N \in \mathbb{R}^{6n} \bigg\vert \exists e_N, \left\lVert\widetilde{\mathcal{C}} \xi_N + \widetilde{\mathcal{D}} e_N\right \rVert_{\Sigma_\nu^{-1}}^2 \leq \delta \right\}, \\
		\Xi_t &= \bigg \{ \xi_t \in \mathbb{R}^{6n} \bigg\vert \exists e_t, \left\lVert\widetilde{\mathcal{C}} \xi_t + \widetilde{\mathcal{D}} e_t\right \rVert_{\Sigma_\nu^{-1}}^2 \leq \delta, \\
		&\qquad \mathcal{A}\xi_t + \mathcal{B} e_t \in \Xi_{t+1} \bigg \}, \: t = 0, \dots, N-1.
	\end{split}
\end{equation}
In the definition of $\Xi_t, t = 0, \dots, N-1$, we have the condition $\mathcal{A} \xi_t + \mathcal{B} e_t \in \Xi_{t+1}$, which ignores the term $\mathcal{K}\nu_{t+1}$. From the structure of $\mathcal{A}$, $\mathcal{C}$, and $\mathcal{K}$, we see that membership in $\Xi_t$ depends only on the $\theta_t$ component of $\xi_t$, which is unaffected by $\mathcal{K} \nu_{t+1}$. That is, we have $\mathcal{A} \xi_t + \mathcal{B} e_t \in \Xi_{t+1}$ if and only if $\mathcal{A} \xi_t + \mathcal{B} e_t + \mathcal{K}\nu_{t+1} \in \Xi_{t+1}$ for any $\nu_{t+1} \in \mathbb{R}^p$. 

We use the sets $\Xi_t$ to determine the existence of recursively feasible attacks at time $t$. 
\begin{lemma}\label{lem: recursiveFeasibility}
	There exists a recursively feasible attack $e_t$ if and only if $\xi_t \in \Xi_t$. That is, there exists a sequence of attacks $\gamma(t, N) = \{e_t, \dots, e_N \}$ such that $\left\lVert \epsilon_t \right \rVert_{\Sigma_\nu^{-1}}^2 \leq \delta, \dots, \left\lVert \epsilon_N \right \rVert_{\Sigma_\nu^{-1}}^2 \leq \delta$ if and only if $\xi_t \in \Xi_t$. 
\end{lemma}

\noindent The proof of Lemma~\ref{lem: recursiveFeasibility} is found in the appendix. The set $\Xi_t$ is nonempty for all $t = 0, \dots, N$ -- {\color{black} if the $\theta_t$ component of $\xi_t$ is equal to 0, then $\xi_t \in \Xi_t$.} This is because, if $\theta_t = 0$, then, following system~\eqref{eqn: eEpsilonSystem}, the attack sequence $\gamma(t, N) = \left\{0, \dots, 0\right \}$ is one such that 
$\left\lVert \epsilon_{t} \right\rVert_{\Sigma_\nu^{-1}}^2 = \dots = \left\lVert \epsilon_{N} \right\rVert_{\Sigma_\nu^{-1}}^2 = 0.$
Recall that system~\eqref{eqn: eEpsilonSystem} has initial state $\theta_0 = 0$, so we have $\xi_0 \in \Xi_0$. This means that the optimization problem~\eqref{eqn: optimalAttack} is feasible for any nonnegative value of $\delta$, i.e., the attacker can always satisfy the detection constraint by choosing not to attack the system.

\section{Attacks Under False Alarm Constraints}\label{sect: equalityConstraints}
In this section, we find an attack sequence $\gamma(0, N)$ that minimizes the cost function $J$ under the constraint that $\left\lVert \epsilon_t \right \rVert_{\Sigma_\nu^{-1}}^2 = 0$, corresponding to finding the optimal attack under the restriction that the probability of being detected is equal to the false alarm probability of the detector. For the case of $\delta = 0$, we can relate the sets $\Xi_t$ to the output minimization problem presented in Lemma~\ref{lem: outputRiccati} and~\cite{Rappaport}. Define
\begin{equation}\label{eqn: gDef}
	\mathcal{G} = \left[\begin{array}{cccc} 0 & I_{3n} & 0 & 0 \end{array} \right].
\end{equation}
The matrix $\mathcal{G}$ selects the variable $\theta_t$ from $\xi_t$ (i.e., $\mathcal{G} \xi_t = \theta_t$). 

\begin{lemma}\label{lem: equalityFeasible}
	{\color{black}For $\delta = 0$ and for $t = 0, \dots, N$, the set $\Xi_t$ is the null space of $\widehat{\mathbf{P}}_{N-t+1}\mathcal{G}$. That is,
	$\Xi_t = \mathscr{N} \left( \widehat{\mathbf{P}}_{N-t+1}\mathcal{G} \right).$}
\end{lemma}
\noindent The proof of Lemma~\ref{lem: equalityFeasible} is found in the appendix. 

The following theorem gives the optimal sequence of attacks when $\delta = 0$. 
\begin{theorem}[Optimal Attack Strategy with $\delta = 0$ Detection Constraint]\label{thm: equalityOptimal}
	An attack sequence $\gamma(0, N)$ that solves~\eqref{eqn: optimalAttack} with $\delta = 0$ is
\begin{align}
	\begin{split}\label{eqn: equalityOptimal2}
	e_t =& -\mathcal{F}_t\left(\mathcal{F}_t^T \mathcal{B}^T \mathbf{Q}_{t+1} \mathcal{B} \mathcal{F}_t \right)^\dagger \mathcal{F}_t^T \mathcal{B}^T\times \\
	& \mathbf{Q}_{t+1} \left(\mathcal{A} - \mathcal{B} \mathcal{D}_t^\dagger \mathcal{C}_t \right) \xi_t - \mathcal{D}_t^\dagger \mathcal{C}_t \xi_t, 
	\end{split}
\end{align}
where
\begin{align}\label{eqn: equalityOptimal3}
	\mathcal{C}_N &= \widetilde{\mathcal{C}}, \quad \mathcal{D}_N = \widetilde{\mathcal{D}}, \quad \mathcal{F}_N = I_s - \mathcal{D}_N^\dagger \mathcal{D}_N,
\end{align}
and, for $t = 0, \dots, N-1$,
\begin{align}
	\begin{split}\label{eqn: equalityOptimal4}
		\mathcal{C}_t &= \left[\begin{array}{c}\widehat{\mathbf{P}}_{N-t} \mathcal{G} \mathcal{A} \\ \widetilde{\mathcal{C}} \end{array} \right], \quad \mathcal{D}_t = \left[\begin{array}{c} \widehat{\mathbf{P}}_{N-t} \widehat{\mathcal{B}} \\ \widetilde{\mathcal{D}} \end{array} \right], 
	\end{split}\\
	\mathcal{F}_t &= I_s - \mathcal{D}_t^\dagger \mathcal{D}_t.
\end{align}
The matrix $\mathbf{Q}_t$ is given recursively backward in time by
\begin{equation}\label{eqn: QtDef}
\begin{split}
		\mathbf{Q}_t &= \mathcal{H}^TQ_t \mathcal{H}+ \left(\mathcal{A} -\mathcal{B} \mathcal{D}_{t}^\dagger \mathcal{C}_t \right)^T \mathbf{Q}_{t+1}\left(\mathcal{A} -\mathcal{B} \mathcal{D}_{t}^\dagger \mathcal{C}_t \right)\\
		& -\left(\mathcal{A} -\mathcal{B} \mathcal{D}_{t}^\dagger \mathcal{C}_t \right)^T \mathbf{Q}_{t+1}\mathcal{B} \mathcal{F}_t \left( \mathcal{F}_t^T \mathcal{B}^T \mathbf{Q}_{t+1} \mathcal{B} \mathcal{F}_t \right)^\dagger \\
		& \times \mathcal{F}_t^T \mathcal{B}^T \mathbf{Q}_{t+1} \left(\mathcal{A} -\mathcal{B} \mathcal{D}_{t}^\dagger \mathcal{C}_t \right),
\end{split}
\end{equation}
with terminal condition
	$\mathbf{Q}_{N+1} = 0$.
\end{theorem}
\noindent Theorem~\ref{thm: equalityOptimal} states that the optimal attack under the $\delta = 0$ detection constraint is a linear feedback of the state $\xi_t$, which is exactly determined by the attacker information $\mathcal{I}_t$. Equation~\eqref{eqn: equalityOptimal2} shows that the optimal attack $e_t$ depends on the matrix $\mathcal{F}_t$, which in turn depends on the matrix $\mathcal{D}_t$. If the matrix $\mathcal{D}_t$ has full column rank, then, $\mathcal{F}_t = 0$, {\color{black} since, by definition, $\mathcal{F}_t$ is the orthogonal projector onto $\mathscr{N} \left( \mathcal{D}_t \right)$.} If the matrix $\mathcal{D}_t$ has full column rank for all $t = 0, \dots, N-1$, then the optimal attack becomes $\gamma(0, N) = \left\{0, \dots, 0 \right\}.$
This corresponds to the case in which the attacker is not powerful enough, and his or her only option to satisfy the $\delta = 0$ detection constraint is to not attack {\color{black} the} system.

Before we prove Theorem~\ref{thm: equalityOptimal}, we provide intermediate results that show that the optimal attack exists and that the optimal attack sequence $e_0, \dots, e_{N-1}$ is unique (the attack $e_N$ may not be unique). The proofs are found in the appendix.

\begin{lemma}\label{lem: QtPSD}
	For all $t = 0, \dots, N$, there exists $\mathcal{U}_t \succeq 0$ such that
		$\mathbf{Q}_t = \mathcal{H}^T Q_t \mathcal{H} + \mathcal{U}_t.$
\end{lemma}

\begin{lemma}\label{lem: solutionExistence}
	For all $t = 0, \dots, N-1$, $\mathscr{R} \left( \mathcal{F}^T_t \mathcal{B}^T\right) = \mathscr{R} \left( \mathcal{F}^T_t \mathcal{B}^T \mathbf{Q}_{t+1} \mathcal{B} \mathcal{F}_t \right)$. 
\end{lemma}
\noindent Define the set
\begin{equation}\label{eqn: zDef}
	Z_t\left(\psi\right) = \left \{ z \in \mathbb{R}^s \vert -\mathcal{F}_t^T \mathcal{B}^T \mathbf{Q}_{t+1} \mathcal{B} \mathcal{F}_t z = \mathcal{F}_t^T \mathcal{B}^T \psi \right \}.
\end{equation}
One consequence of Lemma~\ref{lem: solutionExistence} is that $Z_t \left(\psi\right)$ is nonempty for all $t = 0, \dots, N-1$ and for all $\psi \in \mathbb{R}^{6n}$.

\begin{lemma}\label{lem: zUnique}
	For any $t = 0, \dots, N-1$ and for any $\psi \in \mathbb{R}^{6n}$, if $z_1, z_2 \in Z_t \left(\psi \right)$, then $\mathcal{F}_t z_1 = \mathcal{F}_t z_2$.
\end{lemma}

\begin{IEEEproof}[Proof (Theorem~\ref{thm: equalityOptimal})]
	We resort to dynamic programming to solve~\eqref{eqn: optimalAttack} with $\delta = 0$. The term $\sum_{t=0}^N \trace \left(\widehat{P} Q_t \right)$ in~\eqref{eqn: optimalAttack} does not depend on $\gamma(0, N)$. Define the optimal cost-to-go function for information $\mathcal{I}_t$ as follows:
\begin{align}
	\begin{split}\label{eqn: optimalEqProof2}
		\begin{array}{ccl} J^*_N \left( \mathcal{I}_N \right) =& \underset{e_N}{\text{min}} & \mathbb{E} \left[ \left\lVert\mathcal{H} \xi_N\right\rVert_{Q_N}^2  \Big\vert   \mathcal{I}_N \right] \\ &\text{s.t.} & \epsilon_N = 0 \end{array}
	\end{split},\\
	\begin{split}\label{eqn: optimalEqProof3}
		\begin{array}{ccl} J^*_t \left( \mathcal{I}_t \right) =&\underset{e_t}{\text{min}} & \mathbb{E} \left[ \left\lVert\mathcal{H} \xi_t\right\rVert_{Q_t}^2 + J^*_{t+1} \left( \mathcal{I}_{t+1} \right)  \Big\vert   \mathcal{I}_t \right] \\ &\text{s.t.} & \epsilon_t = 0,\: \mathcal{A} \xi_t + \mathcal{B} e_t \in \Xi_{t+1} \end{array}
	\end{split}.
\end{align}
Equations~\eqref{eqn: optimalEqProof2} and~\eqref{eqn: optimalEqProof3} restrict the attack at each time $t$ to be recursively feasible.

We begin with $t = N$. At time $N$, the attack $e_N$ does not affect the value of $\mathbb{E} \left[\left\lVert\mathcal{H} \xi_N\right\rVert_{Q_N}^2  \right]$, so we choose $e_N$ only to satisfy the constraint $\epsilon_N = \mathcal{C}_N \xi_N + \mathcal{D}_N e_N = 0$. Thus, we have
	$e_N = -\mathcal{D}_N^\dagger \mathcal{C}_N \xi_N,$
and
	$J_N^* \left(\mathcal{I}_N \right) = \xi_N^T\mathbf{Q}_N \xi_N + \Pi_N,$
where
	$\mathbf{Q}_N =  \mathcal{H}^T Q_N \mathcal{H} \text{ and } \Pi_N = 0.$
Proceeding to $N-1$, we first reformulate the constraints. Applying Lemma~\ref{lem: equalityFeasible}, the constraint $\mathcal{A} \xi_{N-1} + \mathcal{B} e_{N-1} \in \Xi_{N}$ becomes
\begin{equation} \label{eqn: optimalEqProof7}
	\widehat{\mathbf{P}}_{1} \mathcal{G} \left( \mathcal{A}\xi_t + \mathcal{B} e_t \right) = 0.
\end{equation}
Combining~\eqref{eqn: optimalEqProof7} with the constraint $\epsilon_{N-1} = 0$ and using the fact that $\mathcal{G}\mathcal{B} = \widehat{\mathcal{B}}$, we have
\begin{equation}\label{eqn: optimalEqProof8}
	\mathcal{C}_{N-1}\xi_{N-1} + \mathcal{D}_{N-1}e_{N-1} = 0,
\end{equation}
where $\mathcal{C}_{N-1}$ and $\mathcal{D}_{N-1}$ are given by~\eqref{eqn: equalityOptimal4}. To solve~\eqref{eqn: optimalEqProof3}, we eliminate the constraint in~\eqref{eqn: optimalEqProof8} (following~\cite{Boyd}) and consider attacks $e_{N-1}$ of the form
\begin{equation}\label{eqn: optimalEqProof10}
	e_{N-1} = \mathcal{F}_{N-1} z_{N-1} - \mathcal{D}_{N-1}^\dagger \mathcal{C}_{N-1} \xi_{N-1},
\end{equation}
where
	$\mathcal{F}_{N-1} = I_s - \mathcal{D}^\dagger_{N-1} \mathcal{D}_{N-1}.$
Equation~\eqref{eqn: optimalEqProof10} describes all recursively feasible $e_{N-1}$ since $\mathscr{R} \left( \mathcal{F}_{N-1} \right) = \mathscr{N} \left( \mathcal{D}_{N-1} \right)$.

After eliminating constraints and performing algebraic manipulations,~\eqref{eqn: optimalEqProof3} becomes
\begin{equation}\label{eqn: optimalEqProof12}
	\begin{split}
	&J^*_{N-1} \left( \mathcal{I}_{N-1} \right) = \xi_{N-1}^T \mathcal{H}^T \mathbf{Q}_{N-1}\mathcal{H} \xi_{N-1} + \Pi_N + \\
	&\quad \trace\left(\Sigma_\nu \mathcal{K}^T \mathbf{Q}_{N} \mathcal{K} \right) +  \min_{z_{N-1}} \overline{\xi}_{N-1}^T \mathbf{Q}_N \overline{\xi}_{N-1},
	\end{split}
\end{equation}
where
	$\overline{\xi}_{N-1} = \left( \mathcal{A} - \mathcal{B} \mathcal{D}^\dagger_{N-1} \mathcal{C}_{N-1} \right)\xi_{N-1} + \mathcal{B} \mathcal{F}_{N-1} z_{N-1}.$
The optimal $z_{N-1}$ satisfies
\begin{equation}\label{eqn: optimalEqProof14}
	0 = \mathcal{F}^T_{N-1} \mathcal{B} \mathbf{Q}_N \overline{\xi}_{N-1}.
\end{equation}
As a consequence of Lemma~\ref{lem: solutionExistence}, such a $z_{N-1}$ exists. One particular $z_{N-1}$ that satisfies~\eqref{eqn: optimalEqProof14} is
\begin{equation}\label{eqn: optimalEqProof15}
\begin{split}
	z_{N-1} =& -\left(\mathcal{F}_{N-1}^T \mathcal{B}^T \mathbf{Q}_N \mathcal{B} \mathcal{F}_{N-1} \right)^\dagger \mathcal{F}_{N-1}^T \mathcal{B}^T \mathbf{Q}_N \times \\
	&\left( \mathcal{A} - \mathcal{B} \mathcal{D}_{N-1}^\dagger \mathcal{C} \right)\xi_{N-1}.
\end{split}
\end{equation}
There may be more than one $z_{N-1}$ that satisfies~\eqref{eqn: optimalEqProof15}. Manipulating~\eqref{eqn: optimalEqProof15}, we have that $z_{N-1}$ satisfies
\begin{equation}\label{eqn: optimalEqProof16}
	-\left(\mathcal{F}_{N-1}^T \mathcal{B}^T \mathbf{Q}_N \mathcal{B} \mathcal{F}_{N-1} \right) z_{N-1} = \mathcal{F}_{N-1}^T \mathcal{B}^T \psi,
\end{equation}
with
$ \psi = \mathbf{Q}_N \left( \mathcal{A} - \mathcal{B} \mathcal{D}_{N-1}^\dagger \mathcal{C} \right)\xi_{N-1}.$
By definition, all $z_{N-1}$ that satisfy~\eqref{eqn: optimalEqProof16} belong to $Z_{N-1} \left( \psi \right)$. Then, since
$e_{N-1} = \mathcal{F}_{N-1} z_{N-1} - \mathcal{D}^\dagger_{N-1} \mathcal{C}_{N-1} \xi_{N-1},$
we have, from Lemma~\ref{lem: zUnique}, that the optimal attack $e_{N-1}$ is unique. 

Substituting~\eqref{eqn: optimalEqProof15} into~\eqref{eqn: optimalEqProof12} and performing algebraic manipulations, we have
\begin{equation}\label{eqn: optimalEqProof17}
	J^*_{N-1} \left( \mathcal{I}_{N-1} \right) = \xi_{N-1}^T \mathcal{H}^T \mathbf{Q}_{N-1} \mathcal{H} \xi_{N-1} + \Pi_{N-1},
\end{equation}
where
\begin{equation}\label{eqn: optimalEqProof18}
	\begin{split}
		&\mathbf{Q}_{N-1} = \mathcal{H}^T Q_{N-1} \mathcal{H} + \left(\mathcal{A} - \mathcal{B} \mathcal{D}_{N-1}^\dagger \mathcal{C}_{N-1} \right)^T \mathbf{Q}_N \times \\
		&\left(\mathcal{A} - \mathcal{B} \mathcal{D}_{N-1}^\dagger \mathcal{C}_{N-1} \right) -  \left(\mathcal{A} - \mathcal{B} \mathcal{D}_{N-1}^\dagger \mathcal{C}_{N-1} \right)^T \times \\
		& \mathbf{Q}_N \mathcal{B} \mathcal{F}_{N-1} \left(\mathcal{F}_{N-1}^T \mathcal{B}^T \mathbf{Q}_N \mathcal{B} \mathcal{F}_{N-1} \right)^\dagger \mathcal{F}_{N-1}^T \mathcal{B}^T \mathbf{Q}_N \times \\
		& \left(\mathcal{A} - \mathcal{B} \mathcal{D}_{N-1}^\dagger \mathcal{C}_{N-1} \right),
	\end{split}
\end{equation}
and
	$\Pi_{N-1} = \Pi_N + \trace\left(\Sigma_\nu \mathcal{K}^T \mathbf{Q}_{N} \mathcal{K} \right).$
Repeating the dynamic programming procedure for $t = N-2, \dots, 0$, we find that the optimal attack has the same form as~\eqref{eqn: optimalEqProof15}, were we replace $N-1$ with $t$.
\end{IEEEproof}

\section{Attacks Under General Detection Constraints}\label{sect: inequalityConstraints}
In this section, we solve~\eqref{eqn: optimalAttack} with positive $\delta$. We design a procedure to find the sequence $\gamma(0, N)$ that minimizes $J$ under the constraint $\left\lVert \epsilon_t \right \rVert_{\Sigma_\nu^{-1}}^2 \leq \delta$ for $t = 0, \dots, N$ (the optimal attack does not have a closed form). This procedure becomes computationally intensive for large $N$. Thus. we also design a less computationally-intensive procedure that finds a sub-optimal and feasible attack sequence.

\subsection{Optimal Attack with $\delta > 0$}\label{sect: ineqOptimal}
For this section only, we introduce the following notation: let $\mathbb{E}_{\xi_t} \left[ \cdot \right]$ denote the expectation taken over $\xi_t$, and let $\mathbb{E}_{\left\{\xi_k \right\}_t^N} \left[ \cdot \right]$ denote the expectation taken over $\xi_t, \dots, \xi_N$. Further, define the operator $\pi_t$ as
	$\pi_t \left( \gamma\left(t, N \right) \right) = e_t.$
That is, $\pi_t$ is an operator that takes an attack sequence over $N-t+1$ time steps and returns the first attack.
To solve~\eqref{eqn: optimalAttack} with $\delta > 0$, we consider, for $t = 0, \dots, N-1$, the problem
\begin{equation}\label{eqn: optimalAttackSubproblem}
	\begin{array}{ccl}\gamma_t^*(t, N) = & \underset{\gamma_t(t, N)}{\text{argmin}} & \mathbb{E}_{\left\{\xi_k \right\}_{t+1}^N} \left[ \sum\limits_{k = t} ^N \left\lVert\mathcal{H} \xi_k\right\rVert_{Q_k}^2 \right] \\ &\text{s.t.} & \left\lVert \epsilon_k \right \rVert_{\Sigma_\nu^{-1}}^2 \leq \delta, k = t, \dots, N
	\end{array},
\end{equation}
where $\gamma_t(t, N) = \left\{e_t, \dots, e_N\right\}$ is an attack sequence in which each attack $e_t, \dots, e_N$ only depends on $\xi_t$\footnote{The state $\xi_t$ is a sufficient statistic for the information set $\mathcal{I}_t$.}. This differs from the definition of $\gamma(t, N)$, in which each attack $e_t, \dots, e_N$ depends on $\xi_t, \dots, \xi_N$, respectively. Problem~\eqref{eqn: optimalAttackSubproblem} has a convex objective and convex constraints, so it can be efficiently solved.

\begin{theorem}[Optimal Attack Strategy with $\delta > 0$ Detection Constraint]\label{thm: inequalityOptimal}
Algorithm~\ref{alg: optimal} gives an attack sequence $\gamma(0, N)$ that solves~\eqref{eqn: optimalAttack} with $\delta > 0$.

\begin{algorithm}[h!]
\caption{Optimal Attack with $\delta > 0$ }\label{alg: optimal}
\begin{algorithmic}[1]
	\State \textbf{Initialize: $\mathcal{I}_0 \gets \left \{ \widetilde{y}_0 \right \}$ }
	\For {$ t = 0, 1, \dots, N-1$}
		\State Solve~\eqref{eqn: optimalAttackSubproblem},  $e_t \gets \pi_t \left(\gamma^*_t \left(t, N\right) \right)$, $\mathcal{I}_{t+1} \gets \left\{ \mathcal{I}_{t}, \widetilde{y}_{t+1}, e_t \right\}$ 
	\EndFor
	\State $e_N \gets \gamma^*_{N-1} \left(N, N \right)$
\end{algorithmic}
\end{algorithm}
\end{theorem}

Algorithm~\ref{alg: optimal} works as follows. At time step $t$, for $t = 0, 1, \dots, N-1$, we find $\gamma_t^* \left(t, N\right)$, the sequence of attacks depending only on $\xi_t$ that solves problem~\eqref{eqn: optimalAttackSubproblem}. The attack $e_t$ is then set as the first attack in the sequence $\gamma_t^*\left(t, N \right)$. In the last ($\left(N+1\right)^{\text{th}}$) time step, the attack $e_N$ is set as the last attack component of the sequence $\gamma_{N-1}^*(N-1, N)$. By construction, every attack $e_t$ produced by Algorithm~\ref{alg: optimal} is recursively feasible: after attacking the system with $e_t$, the subsequence $\gamma_t^*\left(t+1, N\right)$ is a feasible attack sequence at time $t+1$. In order to prove Theorem~\ref{thm: inequalityOptimal}, we require the following Lemma from~\cite{Speyer}:
\begin{lemma}[\!\!\cite{Speyer}]\label{lem: fundamental}
	Let $g(\xi, u)$ be a function such that, for any~$\xi$, $\min_{u \in \mathcal{U}} g(\xi, u)$ exists and $\mathcal{U}$ is a class of functions for which $\mathbb{E}_{\xi} \left[ g(\xi, u) \right]$ exists. Then,
$\min_{u\left(\xi\right) \in \mathcal{U}} \mathbb{E}_\xi \left[ g\left( \xi, u \left( \xi \right) \right) \right] = \mathbb{E}_{\xi} \left[ \min_{u \in \mathcal{U}} g \left(\xi, u \right) \right].$
\end{lemma}

\begin{IEEEproof}[Proof (Theorem~\ref{thm: inequalityOptimal})]
	From problem~\eqref{eqn: optimalAttack}, we have that the optimal cost-to-go function at time $t$, $J^*_t \left(\xi_t \right)$, is defined as
\begin{align}
	\begin{split}\label{eqn: optimalIneqProof1}
		\begin{array}{ccl} J^*_N \left( \xi_N \right) = & \underset{e_N}{\text{min}} &\left\lVert\mathcal{H} \xi_N\right\rVert_{Q_N}^2  \\
		& \text{s.t.} & \left\lVert \epsilon_N \right \rVert_{\Sigma_{\nu}^{-1}}^2 \leq \delta \end{array}
	\end{split},\\
	\begin{split}\label{eqn: optimalIneqProof2}
		\begin{array}{ccl} J^*_t \left( \xi_t \right) &=& \left\lVert\mathcal{H} \xi_t\right\rVert_{Q_t}^2 +  \\
			&\underset{e_t}{\text{min}} & \mathbb{E}_{\xi_{t+1}} \left[  J^*_{t+1} \left( \xi_{t+1} \right)  \vert   \xi_t \right]  \\ &\text{s.t.} & \left\lVert \epsilon_t \right \rVert_{\Sigma_\nu^{-1}}^2 \leq \delta ,\: \mathcal{A} \xi_t + \mathcal{B} e_t \in \Xi_{t+1} \end{array}
	\end{split}.
\end{align}
Let $\widetilde{\gamma}\left(0, N\right) = \left\{ \widetilde{e}_0, \dots, \widetilde{e}_N \right \}$ denote the attack sequence produced by Algorithm~\ref{alg: optimal}. The attack sequence has the form
\begin{equation}\label{eqn: optimalIneqProof3}
	\widetilde{\gamma} \left(0, N \right) = \left \{ \pi_0 \left( \gamma_0^*\left(0, N\right) \right), \dots, \pi_N \left( \gamma_N^*\left(N, N \right) \right) \right\}.
\end{equation}
To show that $\widetilde{\gamma}\left(0, N\right)$ is an optimal attack sequence, we show that each attack $\widetilde{e}_t = \pi_t \left( \gamma^*_t \left( t, N \right) \right)$ is the optimal attack at time $t$, for $t = 0, \dots, N-1$ \footnote{We ignore the attack $\widetilde{e}_N$ because it does not affect the cost associated with $\widetilde{\gamma} \left(0, N\right)$.}.

In order to show that $\widetilde{e}_t$ is the optimal attack, we prove the intermediate result that, for $t = 0, \dots, N-1$,
\begin{equation}\label{eqn: optimalIneqProof4}
\begin{array}{ccl}J_t^*\left(\xi_t \right) = &\underset{\gamma_t(t, N)}{\text{min}} & \mathbb{E}_{\left\{\xi_k \right\}_{t+1}^N} \left[ \sum\limits_{k = t} ^N \left\lVert\mathcal{H} \xi_k\right\rVert_{Q_k}^2 \right] \\ & \text{s.t.} & \left\lVert \epsilon_k \right \rVert_{\Sigma_\nu^{-1}}^2 \leq \delta, k = t, \dots, N
	\end{array}\!\!.
\end{equation}
We resort to induction. In the base case, we show that~\eqref{eqn: optimalIneqProof4} is true for $t=N-1$. 
	($J_N^* \left( \xi_N \right) =\left\lVert\mathcal{H} \xi_N\right\rVert_{Q_N}^2 .$)
Consider the right hand side of~\eqref{eqn: optimalIneqProof4} for $t = N-1$. {\color{black}Expressed in terms of $\xi_{N-1}$,~\eqref{eqn: optimalIneqProof4} becomes}
\begin{align}
	\begin{split}\label{eqn: optimalIneqProof6}
	&\begin{array}{cl} \underset{\gamma_{N-1} \left(N-1, N \right)}{\text{min}} & \mathbb{E}_{\xi_N} \left[\sum
\limits_{k = N-1}^N \left\lVert\mathcal{H} \xi_k\right\rVert_{Q_k}^2  \right] \\ \text{s.t.} & \left\lVert \epsilon_k \right \rVert_{\Sigma_\nu^{-1}}^2 \leq \delta, k = N-1, N\end{array} \\
	&\!\!\!\begin{array}{ccl} = \!\!\!\!& \underset{\gamma_{N-1} \left( N-1, N-1 \right)}{\text{min}} &\!\!\! \left\lVert\mathcal{H} \xi_{N-1}\right\rVert_{Q_{N-1}}^2  \!\! + \mathbb{E}_{\xi_N} \!\!\!\left[ \!\left\lVert\mathcal{H} \xi_N\right\rVert_{Q_N}^2\! \right] \\ & \text{s.t.} & \!\!\!\left\lVert\epsilon_{N-1} \right\rVert_{\Sigma_\nu^{-1}}^2 \leq \delta,  \xi_N \in \Xi_N \end{array},
	\end{split}\\
	\begin{split}\label{eqn: optimalIneqProof7}
		&\!\!\!\begin{array}{ccl} = & \underset{e_{N-1}}{\text{min}} & \left\lVert\mathcal{H} \xi_{N-1}\right\rVert_{Q_{N-1}}^2  +  \mathbb{E}_{\xi_N} \left[ J_N(\xi_N) \vert \xi_{N-1} \right] \\ & \text{s.t.} & \left\lVert\epsilon_{N-1} \right\rVert_{\Sigma_\nu^{-1}}^2 \leq \delta, \\ && \mathcal{A}\xi_{N-1} + \mathcal{B} e_{N-1} \in \Xi_N \end{array},
	\end{split}\\
	&\!\!= J^*_{N-1} \left( \xi_{N-1} \right).
\end{align}
Equation~\eqref{eqn: optimalIneqProof6} follows from the right hand side of~\eqref{eqn: optimalIneqProof4} because the term $ \left\lVert\mathcal{H} \xi_{N-1}\right\rVert_{Q_{N-1}}^2 $ is not affected by $\gamma_{N-1}(N-1, N)$ and because the term $\gamma_{N-1}(N, N)$ does not affect the objective but only needs to satisfy the constraint $\left\lVert\epsilon_{N} \right\rVert_{\Sigma_\nu^{-1}}^2 \leq \delta$. Equation~\eqref{eqn: optimalIneqProof7} follows from~\eqref{eqn: optimalIneqProof6} because the minimization in~\eqref{eqn: optimalIneqProof6} is over $\gamma_{N-1} \left(N-1, N-1 \right)$, which, by definition, is a function of $\xi_{N-1}$. Thus, the expectation over $\xi_N$ in~\eqref{eqn: optimalIneqProof6} refers to the conditional expectation over $\xi_N$ given $\xi_{N-1}$. 

In the induction step, we assume that~\eqref{eqn: optimalIneqProof4} is true for $t+1$ and show that it is true for $t$. The right hand side of~\eqref{eqn: optimalIneqProof4} for $t$ becomes,
\begin{equation}\label{eqn: optimalIneqProof8}
	\begin{array}{cl} \underset{e_t}{\text{min}} & \left[\left\lVert\mathcal{H} \xi_t\right\rVert_{Q_t}^2 + \widetilde{J}_t \left(\xi_t\right)\right] \\	\text{s.t.} & \left\lVert \epsilon_t \right \rVert_{\Sigma_{\nu}^{-1}}^2 \leq \delta, \mathcal{A}\xi_t + \mathcal{B} e_t \in \Xi_{t+1} \end{array},
\end{equation}
where  $\widetilde{J}_t \left(\xi_t\right)$ is defined as
\begin{equation}\label{eqn: optimalIneqProof9}
	\begin{array}{cl} \underset{\gamma_t(t+1, N)}{\text{min}} & \mathbb{E}_{\left\{\xi_k \right\}_{t+1}^N} \left[ \sum\limits_{k = t+1} ^N\left\lVert\mathcal{H} \xi_k\right\rVert_{Q_k}^2 \right] \\ \text{s.t.} & \left\lVert \epsilon_k \right \rVert_{\Sigma_\nu^{-1}}^2 \leq \delta, k = t+1, \dots, N\end{array}.
\end{equation}
The expression in~\eqref{eqn: optimalIneqProof8} follows from the right hand side of~\eqref{eqn: optimalIneqProof4} because $\gamma_t(t, N)$ can be partitioned as $\left\{ e_t, \gamma_t \left( t+1, N \right) \right\}$, $\xi_t$ does not depend on $\gamma_t(t+1, N)$, and a feasible $\gamma_t \left(t+1, N\right)$ in~\eqref{eqn: optimalIneqProof9} exists if and only if $e_t$ is recursively feasible. 

Further manipulating~\eqref{eqn: optimalIneqProof9}, we have
\begin{align}
\begin{split}\label{eqn: optimalIneqProof10}
	&\widetilde{J}_t\left(\xi_t\right) = \mathbb{E}_{\xi_{t+1}} \left[\widehat{J}_{t+1} \left(\xi_{t+1} \right) \Big \vert \xi_{t} \right],
\end{split}
\end{align}
where $\widehat{J}_{t+1} \left(\xi_{t+1} \right)$ is defined as
\begin{equation}\label{eqn: optimalIneqProof11}
\begin{split}
	\begin{array}{cl} \underset{\gamma_{t+1} \left(t+1, N \right)}{\text{min}} & \mathbb{E}_{\left\{\xi_k \right\}_{t+2}^N}\left[ \sum\limits_{k = t+1}^N \left\lVert\mathcal{H} \xi_k\right\rVert_{Q_k}^2 \right] \\ \text{s.t.} & \left\lVert \epsilon_k \right \rVert_{\Sigma_\nu^{-1}}^2 \leq \delta, k = t+1, \dots, N\end{array}.
\end{split}
\end{equation}
Equation~\eqref{eqn: optimalIneqProof10} follows from~\eqref{eqn: optimalIneqProof9} because the minimization in~\eqref{eqn: optimalIneqProof9} is over $\gamma_t \left(t+1, N\right)$, which, by definition, is a function of $\xi_t$, so the expectation over $\xi_{t+1}, \dots, \xi_{N}$ refers to the conditional expectation given $\xi_t$. We also use Lemma~\ref{lem: fundamental} to exchange the minimization operation and expectation over $\xi_{t+1}$ in~\eqref{eqn: optimalIneqProof9} to derive~\eqref{eqn: optimalIneqProof10}. By the induction hypothesis, we have
$\widehat{J}_{t+1} \left(\xi_{t+1}\right) = J_{t+1}^* \left(\xi_{t+1} \right).$
Substituting back into~\eqref{eqn: optimalIneqProof10} and~\eqref{eqn: optimalIneqProof8} shows that~\eqref{eqn: optimalIneqProof4} is true for $t$. 

To conclude the proof of Theorem~\ref{thm: inequalityOptimal}, we note that, as a result of~\eqref{eqn: optimalIneqProof4}, for $t = 0, \dots, N-1$,
\begin{equation}\label{eqn: optimalIneqProof12}
	\begin{array}{ccl} \widetilde{e}_t = &\underset{e_t}{\text{argmin}} & \mathbb{E}_{\xi_{t+1}} \left[  J^*_{t+1} \left( \xi_{t+1} \right)  \vert   \xi_t \right]  \\ &\text{s.t.} & \left\lVert \epsilon_t \right \rVert_{\Sigma_\nu^{-1}}^2 \leq \delta ,\: \mathcal{A} \xi_t + \mathcal{B} e_t \in \Xi_{t+1} \end{array}.
\end{equation}
Because $\widetilde{e}_t$ is optimal for $t = 0, \dots, t=N-1$ and $e_N$ does not affect the cost, $\widetilde{\gamma} \left( 0, N \right)$ is an optimal attack sequence.
\end{IEEEproof}

\subsection{Windowed Attack Algorithm with $\delta > 0$}\label{sect: windowedIneq}
Although Algorithm~\ref{alg: optimal} gives the optimal attack sequence with $\delta>0$, it is also computationally expensive. In the first time step ($t = 0$), in order to find the optimal attack, Algorithm~\ref{alg: optimal} aims to find an optimal $N$-length attack sequence subject to $N$ constraints. In each subsequent time step, the length of the attack sequence over which the optimization occurs and the number of constraints in the optimization only decreases by~$1$. Even though~\eqref{eqn: optimalAttackSubproblem} is a convex optimization, if $N$ is large, then, in order to find the optimal attack, Algorithm~\ref{alg: optimal} must repeatedly solve large (convex) optimization problems, each with a large number of constraints.

To find an attack sequence in a less computationally-intensive manner, we consider the windowed attacker optimization problem:
\begin{equation}\label{eqn: windowedOptimizationProblem}
\begin{split}
	&\widehat{\gamma}_t(t, t+W-1) = \\
	&\begin{array}{cl}\underset{\gamma_t(t, t+W-1)}{\text{argmin}} & \mathbb{E}_{\left\{\xi_k \right\}_{t+1}^{t+W-1}} \left[ \sum\limits_{k = t} ^{t+W-1}\left\lVert\mathcal{H} \xi_k\right\rVert_{Q_k}^2 \right] \\ \text{s.t.} &  \left\lVert \epsilon_k \right \rVert_{\Sigma_\nu^{-1}}^2 \leq \delta, k = t, \dots, t+W-1, \\
& \left\lVert \widehat{\mathbf{P}}_{N-(t+W-1)} \mathcal{G} \xi_{t+W} \right\rVert_{\Sigma_\nu^{-1}}^2 \leq \delta, \end{array}
\end{split} 
\end{equation}
where $W \in \left\{2, \dots, N+1\right\}$ is a predetermined window size chosen by the attacker\footnote{We do not consider the case of $W = 1$ since, for any time $t$, the current attack $e_t$ does not affect the value of $\xi_t$.}. The goal of problem~\eqref{eqn: windowedOptimizationProblem} is to find a $W$-length attack sequence that minimizes the attacker's cost over $W$ time steps. 

We use problem~\eqref{eqn: windowedOptimizationProblem} to find a suboptimal attack sequence in less computationally intensive manner than Algorithm~\ref{alg: optimal}.

\begin{algorithm}[h!]
\caption{Windowed Attack with $\delta > 0$ }\label{alg: windowed}
\begin{algorithmic}[1]
	\State \textbf{Initialize: $\mathcal{I}_0 \gets \left \{ \widetilde{y}_0 \right \}$ }
	\For {$t = 0, 1, \dots, N-W$}
		\State Solve~\eqref{eqn: windowedOptimizationProblem}, $e_t \gets \pi_{N-W+1} \left( \widehat{\gamma} \left(t, t+W-1 \right) \right)$,  $\mathcal{I}_{t+1} \gets \left\{ \mathcal{I}_{t}, \widetilde{y}_{t+1}, e_t \right\}$ 
	\EndFor
	\For {$ t = N-W+1, \dots, N$}
		\State Solve~\eqref{eqn: optimalAttackSubproblem}, $e_t \gets \pi_t \left(\gamma^*_t \left(t, N\right) \right)$, $\mathcal{I}_{t+1} \gets \left\{ \mathcal{I}_{t}, \widetilde{y}_{t+1}, e_t \right\}$ 
	\EndFor
	\State $e_N \gets \gamma^*_{N-1} \left(N, N \right)$
\end{algorithmic}
\end{algorithm}
\noindent Algorithm~\ref{alg: windowed} works as follows. At each time step $t$, for $t = 0, \dots, N-W$, we find a $W$-length attack sequence (that depends only on $\xi_t$) that minimizes the attacker's cost over $W$ time steps. We set the attack $e_t$ to be the first component of the attack sequence $\widehat{\gamma}_t \left( t, t+W-1 \right)$. At each time step $t$, for $t = N-W+1, \dots, N$, we determine the attack $e_t$ in the same way as in Algorithm~\ref{alg: optimal}.

The additional constraint 
	$\left\lVert \widehat{\mathbf{P}}_{N-(t+W-1)} \mathcal{G} \xi_{t+W} \right\rVert_{\Sigma_\nu^{-1}}^2 \leq \delta$
ensures that the attack $e_t$, as determined by Algorithm~\ref{alg: windowed}, is recursively feasible. The attack $e_t$ is the first component of a sequence $\widehat{\gamma}_t \left(t, t+W-1\right)$ that satisfies all constraints in~\eqref{eqn: windowedOptimizationProblem}. The additional constraint requires that, for the state $\xi_{t+W}$, which depends on $\widehat{\gamma}_t \left(t, t+W-1\right)$, there exists an attack sequence $\widehat{\gamma}_t \left(t+W, N \right)$ such that  $\sum_{k = t+W}^N \left\lVert \epsilon_k \right \rVert_{\Sigma_\nu^{-1}}^2 \leq \delta$. Since the sum of $\left\lVert \epsilon_k \right \rVert_{\Sigma_\nu^{-1}}^2 $ from $k = t+W$ to $k = N$ is no greater than $\delta$, we have $\left\lVert \epsilon_k \right \rVert_{\Sigma_\nu^{-1}}^2  \leq \delta$ for all $k = t+W, \dots, N$. Thus, the attack sequence
$\left\{\widehat{\gamma}_t \left(t, t+W-1 \right), \widehat{\gamma}_t \left(t+W, N \right) \right\}$
satisfies  $\left\lVert \epsilon_k \right \rVert_{\Sigma_\nu^{-1}}^2  \leq \delta$ for all $k = t, \dots, N$, so the attack $e_t$ must be recursively feasible. 

In each time step of Algorithm~\ref{alg: windowed}, we minimize cost over an attack sequence of length no greater than $W$, and subject to no more than $W+1$ constraints. Thus, even if $N$ is large, an attacker can use Algorithm~\ref{alg: windowed} to find a feasible attack sequence without the computational expense of Algorithm~\ref{alg: optimal} by choosing $W$ to be small. One drawback of a small window size $W$ is that the resulting attack sequence incurs a larger cost $J$. We verify the trade off between window size and optimality gap via numerical simulation in Section~\ref{sect: examples}. 

\section{Numerical Examples}\label{sect: examples}
We demonstrate the proposed attack strategies under detection constraints with separate examples for the $\delta = 0$ and $\delta > 0$ cases. We consider the linearized state space model of a helicopter provided by~\cite{Helicopter}. The model is comprised of $10$ states and $4$ actuators. Due to space constraints, we refer the reader to~\cite{Helicopter} for a detailed explanation of the model states and the numerical values of the $A$ and $B$ matrices. The helicopter has sensors for each of the state variables ($C = I_{10}$). In our examples, we consider the following statistical properties: $\overline{x}_{t_0} = 0, \Sigma_x = 5 I_{10}, \Sigma_v = 10^{-3} I_{10}, \Sigma_w = 10^{-4} \text{diag} \left(6, .1, 2, 2, .1, 6, 2, 2, 2, .1 \right).$ 

For all numerical examples, the attacker has target state
\[ x^* = \left[\begin{array}{cccccccccc} 0 & 4 & 0 & 0 & 0 & 8.2 & 0 & 0 & 0 & 0 \end{array}\right]^T.\]
The attacker has the following $Q_t$ matrix for all $t$:
\[Q_t = \diag (\begin{array}{cccccccccc} .1,& 3,& .1,& .1,& .1, &4, &.1, &.1, &.1, &.1\end{array}), \]
which means that the attacker only cares about manipulating the $x_{(2)}$ and $x_{(6)}$ components of the helicopter's state, corresponding to vertical and lateral velocity, respectively. The system starts running at $t = -75$, and the attacker attacks the system from $t = 0$ to $t = 75$ (i.e., $N = 75$). 

\begin{figure}[h!]
	\centering
	\includegraphics[keepaspectratio = true, scale = .6]{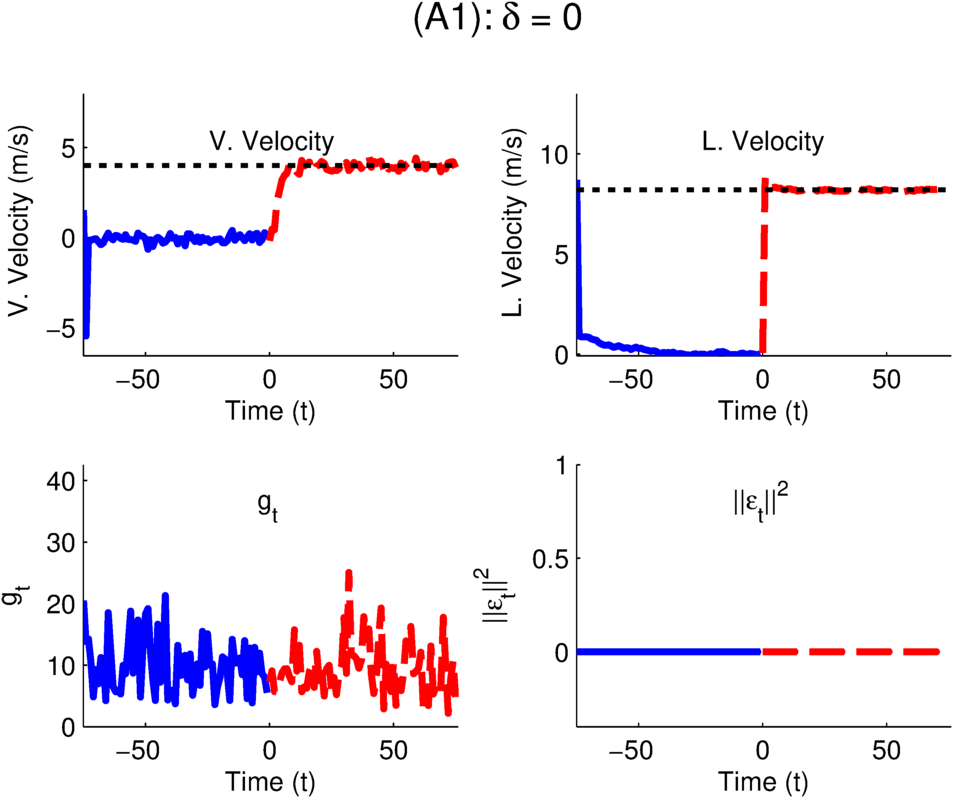}
	\caption{Effect of the optimal attack (A1) under $\delta = 0$ constraint.  Top: system states versus time. The black dotted line is the target state. Bottom: detection statistic versus time}
	\label{fig: numericalResults1}
\end{figure}

\begin{figure}[h!]
	\centering
	\includegraphics[keepaspectratio = true, scale = .6]{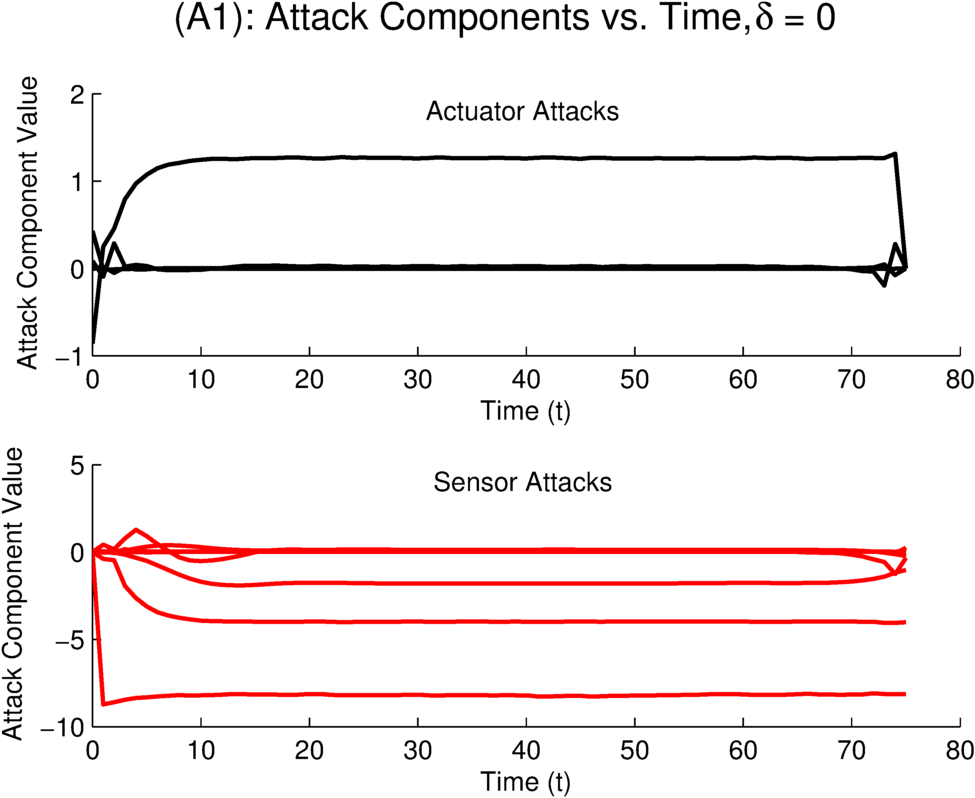}
	\caption{\color{black} Component-wise description of the optimal attack (A1) under the $\delta = 0$ constraint.}
	\label{fig: A1Attack}
\end{figure}

First, we consider an attacker, denoted as ``A1'', that can attack all of the actuators and eight of the sensors -- the attacker cannot alter the sensors measuring $x_{(8)}$ (yaw rate) and $x_{(9)}$ (roll angle). Figure~\ref{fig: numericalResults1} shows the effect of the optimal attack (A1) with the $\delta = 0$ constraint, {\color{black} and Figure~\ref{fig: A1Attack} provides a component-wise description of the optimal attack over time.}
From time $t = 0$ to $t = 75$, the attack computed using Theorem~\ref{thm: equalityOptimal} moves the system to the target state while satisfying the $\left\lVert \epsilon_t \right \rVert_{\Sigma_{\nu}^{-1}}^2 = 0$ constraint. 

Second, we consider an attacker, denoted as ``A2'', that can attack inputs $u_{(3)}$ and $u_{(4)}$ and manipulate the sensor values measuring $x_{(2)}, x_{(6)}, x_{(7)}, \text{ and } x_{(10)}$. 
\begin{figure}[h!]
	\centering
	\includegraphics[keepaspectratio = true, scale = .6]{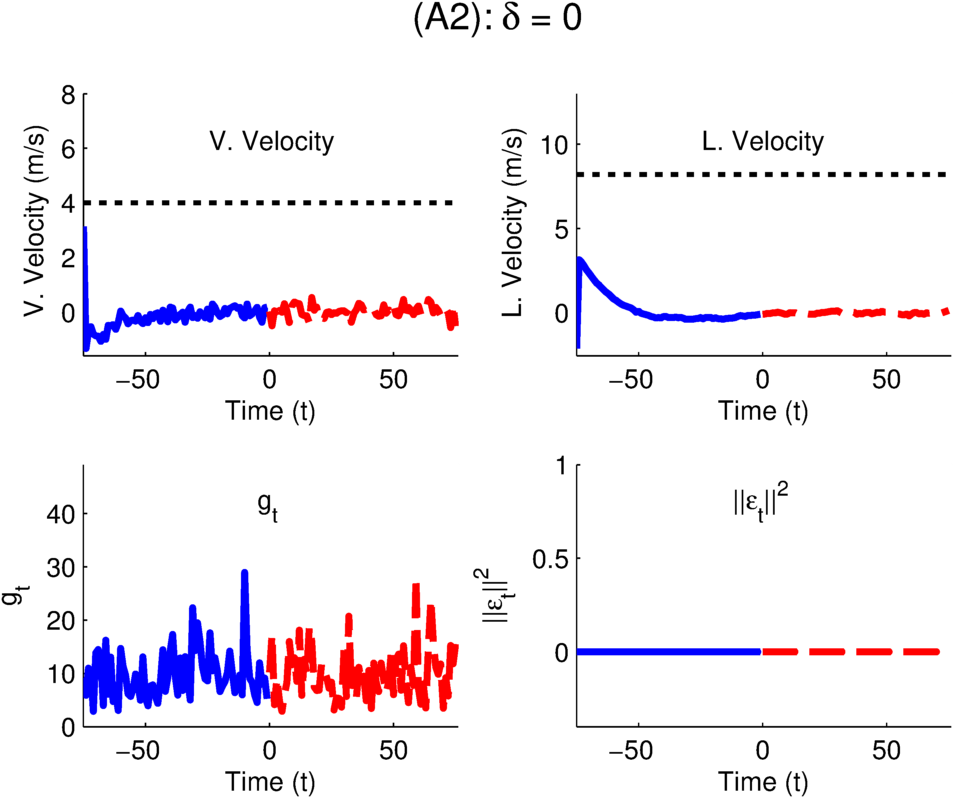}
	\caption{Effect of the optimal attack (A2) under $\delta = 0$ constraint.  Top: system states versus time. The black dotted line is the target state. Bottom: detection statistic versus time}
	\label{fig: numericalResults3}
\end{figure}
\begin{figure}[h!]
	\centering
	\includegraphics[keepaspectratio = true, scale = .6]{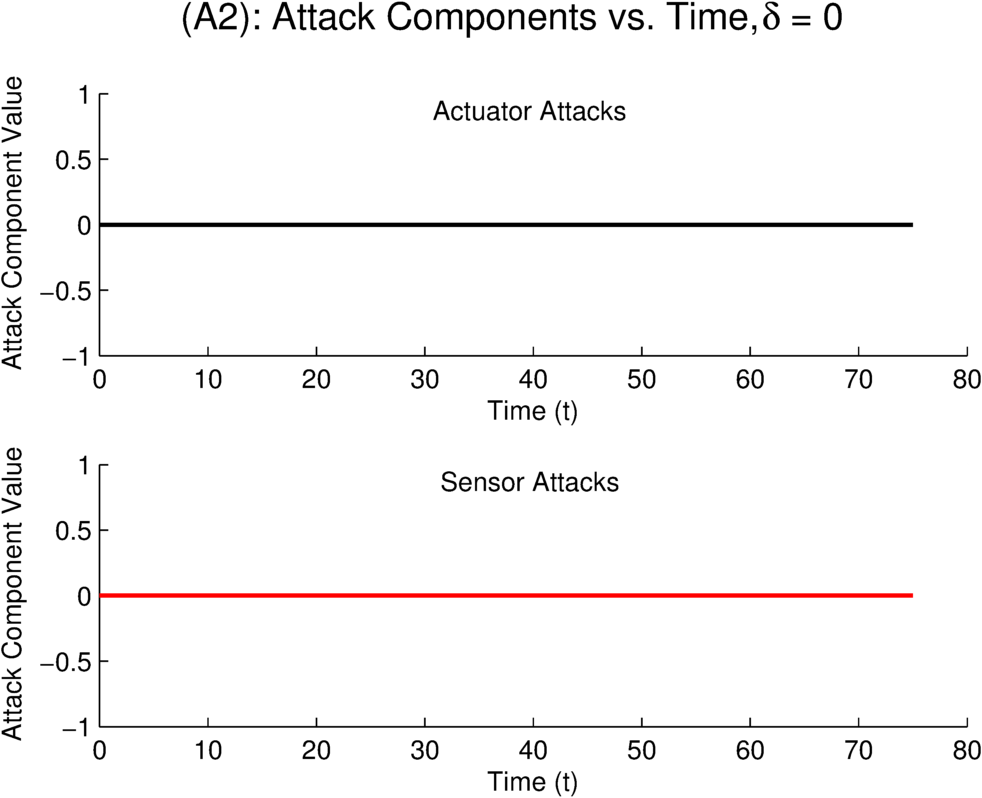}
	\caption{\color{black} Component-wise description of the optimal attack (A2) under the $\delta = 0$ constraint. The only feasible attack for (A2) is to not attack the system.}
	\label{fig: A2Attack}
\end{figure}
\noindent Figure~\ref{fig: numericalResults3} shows that, for the attacker (A2), the optimal attack with $\delta = 0$ constraint does not successfully move the system to the target state. This is because the attacker is not powerful enough and cannot attack enough sensors and actuators. {\color{black} Indeed, as Figure~\ref{fig: A2Attack} shows, the optimal strategy for (A2) under the $\delta = 0$ constraint is to not attack the system.}

For demonstrating the attack strategies with $\delta > 0$, we consider the attacker (A2). In the implementation of Algorithms~\ref{alg: optimal} and~\ref{alg: windowed}, we solve the optimization problems~\eqref{eqn: optimalAttackSubproblem} and~\eqref{eqn: windowedOptimizationProblem} using MOSEK~\cite{mosek}. First, we consider the optimal attack from Algorithm~\ref{alg: optimal} under the constraint $\delta = 1$.
\begin{figure}[h!]
	\centering
	\includegraphics[keepaspectratio = true, scale = .6]{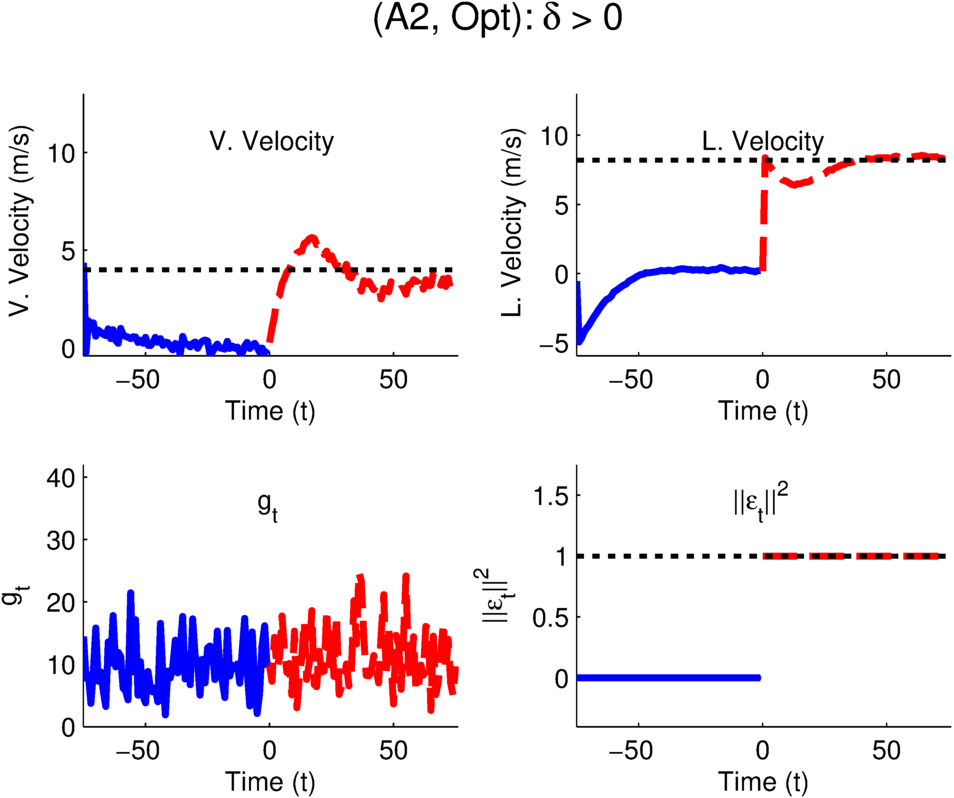}
	\caption{Effect of the optimal attack (A2) under $\delta > 0$ constraint.  Top: system states versus time. The black dotted line is the target state. Bottom: detection statistic versus time. The black dotted line is $\delta$.}
	\label{fig: numericalResults5}
\end{figure}
\begin{figure}[h!]
	\centering
	\includegraphics[keepaspectratio = true, scale = .6]{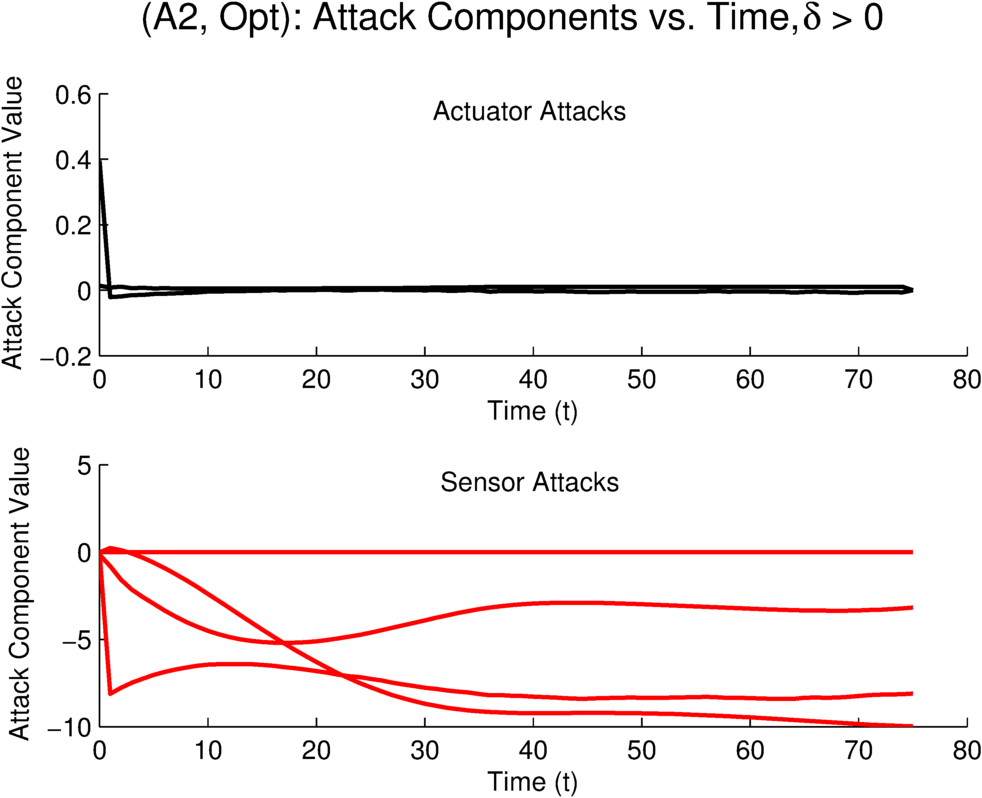}
	\caption{\color{black} Component-wise description of the optimal attack (A2) under the $\delta > 0$ constraint.}
	\label{fig: IneqOptA2Attack}
\end{figure}
\noindent Figure~\ref{fig: numericalResults5} shows that, by following Algorithm~\ref{alg: optimal}, the attacker (A2) is able to move the system to the target state while satisfying the $\left \lVert \epsilon_t \right\rVert_{\Sigma_\nu^{-1}}^2 \leq 1$ constraint. {\color{black} Figure~\ref{fig: IneqOptA2Attack} provides a component-wise description of the optimal attack under the inequality constraint.} 

\begin{figure}[h!]
	\centering
	\includegraphics[keepaspectratio = true, scale = .6]{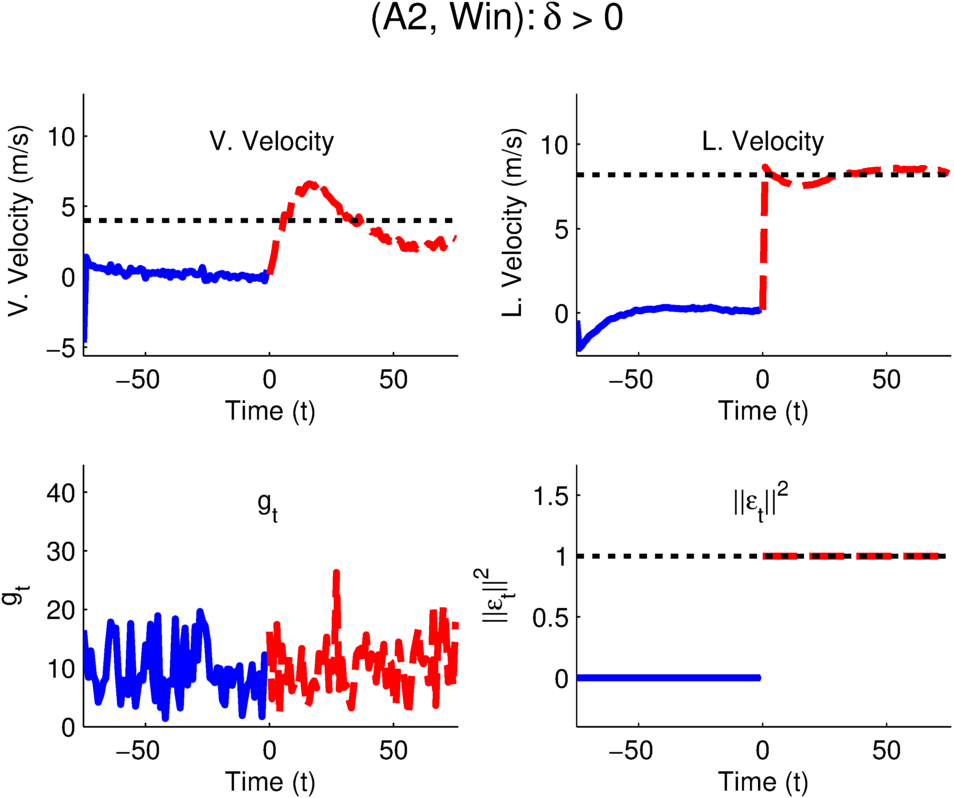}
	\caption{Effect of the windowed attack (A2) under $\delta > 0$ constraint. Top: system states versus time. The black dotted line is the target state. Bottom: detection statistic versus time. The black dotted line is $\delta$.}
	\label{fig: numericalResults7}
\end{figure}

\begin{figure}[h!]
	\centering
	\includegraphics[keepaspectratio = true, scale = .6]{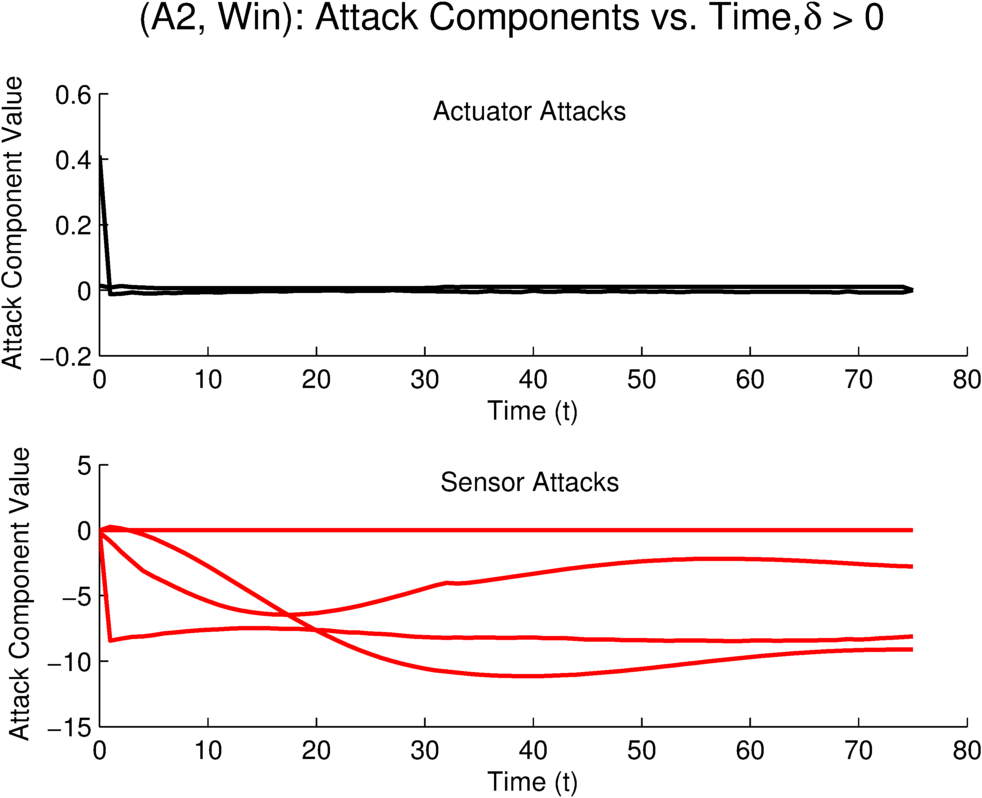}
	\caption{\color{black} Component-wise description of the windowed attack (A2) under the $\delta > 0$ constraint.}
	\label{fig: IneqWindowedA2Attack}
\end{figure}

We then consider the attacker (A2) using the windowed attack algorithm (Algorithm~\ref{alg: windowed}) with window size $W = 5$. Figure~\ref{fig: numericalResults7} shows that, like the optimal attack, the suboptimal attack computed using Algorithm~\ref{alg: windowed} also successfully brings the system state to the target state while satisfying the $\left \lVert \epsilon_t \right\rVert_{\Sigma_\nu^{-1}}^2 \leq 1$ constraint. Moreover, {\color{black} Figure~\ref{fig: IneqWindowedA2Attack} shows that the windowed attack resembles the optimal inequality constrained attack. Both the optimal and windowed attack satisfy the $\left\lVert \epsilon_t \right \rVert_{\Sigma_{\nu}^{-1}}^2 \leq 1$ constraint with inequality, which shows that the optimal cost is achieved at the boundary of the constraint set. The optimal attack strategy induces as large a bias $\epsilon_t$ as allowed by the explicit detection avoidance constraints.}

The performance of Algorithm~\ref{alg: windowed} depends on the window size $W$. We evaluate the cost attained by Algorithm~\ref{alg: windowed} (by attacker (A2)) as a function of the window size $W$ and the value of the constraint bound $\delta$. For each $W$ and each value of $\delta$, we compute the optimality gap of Algorithm~\ref{alg: windowed} as the average cost of $10000$ simulations. The optimal cost is the one obtained by Algorithm~\ref{alg: optimal}, which we compute as the average of $10000$ simulations.
\begin{figure}[h!]
	\centering
	\includegraphics[keepaspectratio = true, scale = .6]{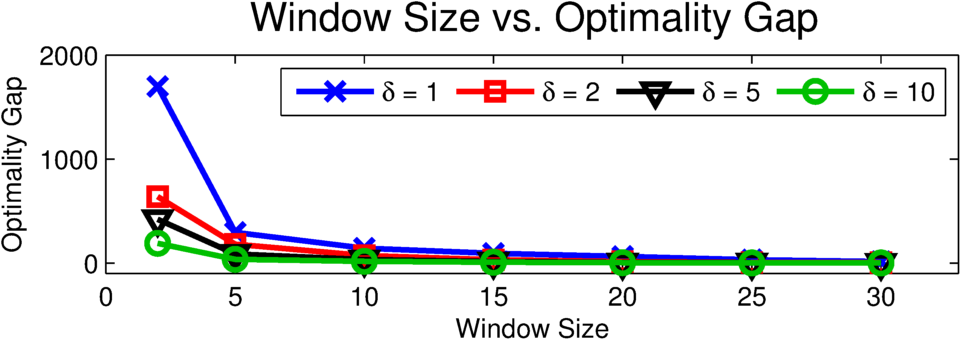}
	\caption{Performance of Algorithm~\ref{alg: windowed} as a function of window size.}
	\label{fig: numericalResults9}
\end{figure}
{\color{black}In general, for a fixed target state $x^*$, the optimality gap between the optimal attack and windowed attack depends on the window size $W$ and the bound $\delta$.} Figure~\ref{fig: numericalResults9} shows that, as window size $W$ increases, the cost achieved by Algorithm~\ref{alg: windowed} approaches the optimal cost. {\color{black} For the smallest window size $W = 2$, which corresponds to the greedy attack strategy of minimizing the one-step-ahead cost under recursive feasibility constraints, tighter detection avoidance constraints (i.e., lower values of $\delta$) incur a noticeably larger optimality gap. One possible explanation is that, for tighter constraints, there is a limited set of feasible attacks at any time step, and the greedy strategy further limits the set of feasible attacks at future time-steps, resulting in a larger overall optimality gap. That is, the impact of the greedy strategy on the optimality gap is greater for tighter constraints.}

\section{Conclusion}\label{sect: conclusion}
In this paper we studied attackers with control objectives and detection constraints against CPS. We formulated a cost function that captures the attacker's control objectives and defined constraints that relate to the probability of the attack being detected. In the case that the attacker's probability of being detected is constrained to be the false alarm rate of the detector, we showed that the optimal attack strategy is a linear feedback of an augmented system state calculated from the attacker's information set. Under more general constraints to the attacker's effect on the CPS's detection statistic, we provided an algorithm to find the optimal attack sequence and a second, less computationally intensive, algorithm to find a feasible, sub-optimal attack. Finally, we illustrated our attack strategies through numerical examples involving a remotely-controlled helicopter under attack. 

\appendix
\subsection{Proof of Lemma~\ref{lem: recursiveFeasibility}}\label{sect: recursiveFeasibilityProof}
\begin{IEEEproof}
	We resort to induction. The base case of $t = N$ is true by the definition of $\Xi_N$. In the induction step, we assume there exists a sequence of attacks $\gamma(t+1, N)$ such that $\left\lVert \epsilon_{t+1} \right \rVert_{\Sigma_\nu^{-1}}^2 \leq \delta, \dots, \left\lVert \epsilon_N \right \rVert_{\Sigma_\nu^{-1}}^2 \leq \delta$ if and only if $\xi_{t+1} \in \Xi_{t+1}$, and we show that there exists a recursively feasible $e_t$ if and only if $\xi_t \in \Xi_{t}$. 

\textit{(If)} Let $\xi_t \in \Xi_{t}$. Then there exists $e_t$ such that $\left\lVert\widetilde{\mathcal{C}} \xi_t + \widetilde{\mathcal{D}} e_t\right \rVert_{\Sigma_\nu^{-1}}^2 \leq \delta$ and $\xi_{t+1} = \mathcal{A} \xi_t + \mathcal{B} e_t + \mathcal{K} \nu_{t+1} \in \Xi_{t+1}$. Since $\xi_{t+1} \in \Xi_{t+1}$, by the induction hypothesis, there exists $\gamma(t+1, N)$ such that $\left\lVert \epsilon_{t+1} \right \rVert_{\Sigma_\nu^{-1}}^2 \leq \delta, \dots, \left\lVert \epsilon_N \right \rVert_{\Sigma_\nu^{-1}}^2 \leq \delta$. Concatenating $e_t$ and $\gamma(t+1, N)$, we have that $\gamma(t, N) = \{ e_t, \gamma(t+1, N)\}$ is an attack such that $\left\lVert \epsilon_{t} \right \rVert_{\Sigma_\nu^{-1}}^2 \leq \delta, \dots, \left\lVert \epsilon_N \right \rVert_{\Sigma_\nu^{-1}}^2 \leq \delta$, which means that $e_t$ is recursively feasible.

\textit{(Only If)} Let $e_t$ be a recursively feasible attack. Then, there exists $\gamma(t, N) = \{ e_t, \gamma(t+1, N)\}$ such that $\left\lVert \epsilon_{t} \right \rVert_{\Sigma_\nu^{-1}}^2 \leq \delta, \dots, \left\lVert \epsilon_N \right \rVert_{\Sigma_\nu^{-1}}^2 \leq \delta$. Since the subsequence $\gamma(t+1, N)$ satisfies $\left\lVert \epsilon_{t+1} \right \rVert_{\Sigma_\nu^{-1}}^2 \leq \delta, \dots, \left\lVert \epsilon_N \right \rVert_{\Sigma_\nu^{-1}}^2 \leq \delta$, we have, by the induction hypothesis, that $\xi_{t+1} = \mathcal{A} \xi_t + \mathcal{B}e_t + \mathcal{K}\nu_{t+1} \in \Xi_{t+1}$. Since $\left\lVert \epsilon_t \right\rVert_{\Sigma_\nu^{-1}}^2 \leq \delta$, we have $\left\lVert\widetilde{\mathcal{C}} \xi_t + \widetilde{\mathcal{D}} e_t\right \rVert_{\Sigma_\nu^{-1}}^2 \leq \delta$. This means that $\xi_t \in \Xi_t$. 
\end{IEEEproof}

\subsection{Proof of Lemma~\ref{lem: equalityFeasible}}\label{sect: equalityFeasibleProof}
\begin{IEEEproof}
	We resort to induction. In the base case, we show that $\Xi_N = \mathscr{N} \left( \widehat{\mathbf{P}}_1 \mathcal{G} \right)$. By definition of $\Xi_N$, we have $\xi_N \in \Xi_N$ if and only if there exists $e_N$ such that $\left\lVert\widetilde{\mathcal{C}} \xi_N + \widetilde{\mathcal{D}} e_N\right \rVert_{\Sigma_\nu^{-1}}^2 =0$. 
Applying the result of Lemma~\ref{lem: outputRiccati} and noting that $\mathcal{G}\xi_N = \theta_N$, we have that $ \left\lVert \widehat{\mathcal{C}}\mathcal{G} \xi_N + \widetilde{\mathcal{D}} e_N\right \rVert_{\Sigma_\nu^{-1}}^2 =0$ if and only if $\widehat{\mathbf{P}}_1\mathcal{G} \xi_N = 0$. Thus, we have $\xi_N \in \Xi_N$ if and only if $\widehat{\mathbf{P}}_1\mathcal{G} \xi_N = 0$, which shows that $\Xi_N = \mathscr{N} \left( \widehat{\mathbf{P}}_1 \mathcal{G} \right).$

In the induction step, we assume that $\Xi_{t+1} = \mathscr{N} \left( \widehat{\mathbf{P}}_{N-t} \mathcal{G}\right),$ and we show that $\Xi_{t} = \mathscr{N} \left( \widehat{\mathbf{P}}_{N-t+1} \mathcal{G} \right).$ By definition of $\Xi_t$, we have $\xi_t \in \Xi_t$ if and only if there exists $e_t$ such that $\left\lVert\widetilde{\mathcal{C}} \xi_t + \widetilde{\mathcal{D}} e_t\right \rVert_{\Sigma_\nu^{-1}}^2 =0$ and $\xi_{t+1} = \mathcal{A}\xi_t + \mathcal{B} e_t + \mathcal{K} \nu_{t+1} \in \Xi_{t+1}$. By the induction hypothesis $\Xi_{t+1} = \mathscr{N} \left( \widehat{\mathbf{P}}_{N-t} \mathcal{G} \right)$. Thus, we have $\xi_{t+1} \in \Xi_{t+1}$ if and only if $\widehat{\mathbf{P}}_{N-t} \mathcal{G}\xi_{t+1} = 0$. Applying the results of {\color{black}Lemma~\ref{lem: recursiveFeasibility}}, we then have that $\xi_{t+1} \in \Xi_{t+1}$ if and only if there exists an attack sequence $\gamma(t+1, N)$ such that, starting from state $\xi_{t+1}$, we have
$\left\lVert \epsilon_{t+1} \right\rVert_{\Sigma_\nu^{-1}}^2 = \dots = \left\lVert \epsilon_{N} \right\rVert_{\Sigma_\nu^{-1}}^2 = 0.$
Concatenating $e_t$ and $\gamma(t+1, N)$, we have that $\gamma(t, N) = \left\{e_t, \gamma(t+1, N) \right \}$ is an attack sequence such that, starting from state $\xi_t$, we have
$\left\lVert \epsilon_{t} \right\rVert_{\Sigma_\nu^{-1}}^2 = \dots = \left\lVert \epsilon_{N} \right\rVert_{\Sigma_\nu^{-1}}^2 = 0.$
Since there exists such an attack sequence $\gamma(t, N)$, we have
	$\min_{\gamma(t, N)} \sum_{k=t}^N \left\lVert \epsilon_k \right\rVert_{\Sigma_\nu^{-1}}^2 = 0.$
Applying the results of Lemma~\ref{lem: outputRiccati}, we have $\xi_t \in \Xi_t$ if and only if $\widehat{\mathbf{P}}_{N-t+1}\mathcal{G}\xi_t = 0$, which shows that $\Xi_{t} = \mathscr{N} \left( \widehat{\mathbf{P}}_{N-t+1} \mathcal{G} \right).$
\end{IEEEproof}

\subsection{Proof of Lemma~\ref{lem: QtPSD}}\label{sect: QtPSDProof}
\begin{IEEEproof}
	We resort to induction. In the base case $t = N$, we have $ \mathbf{Q}_N = \mathcal{H}^T Q_N \mathcal{H}$ by definition. In the induction step, we assume that 
		$\mathbf{Q}_{t+1} = \mathcal{H}^T Q_{t+1} \mathcal{H} + \mathcal{U}_{t+1}$
	for some $\mathcal{U}_{t+1} \succeq 0$, and we show that there exists $\mathcal{U}_{t} \succeq 0$ such that
		$\mathbf{Q}_{t} = \mathcal{H}^T Q_t \mathcal{H} + \mathcal{U}_t.$
From the induction hypothesis, we have that $\mathbf{Q}_{t+1} \succeq 0$ since ${\color{black}Q_{t+1}} \succ 0$ and $\mathcal{U}_{t+1} \succeq 0$. Then, by algebraic manipulation of equation~\eqref{eqn: QtDef}, we have
\begin{equation}\label{eqn: psdProof3}
	\mathbf{Q}_t = \mathcal{H}^T Q_t \mathcal{H} + \mathcal{X}_t^T \mathbf{Q}_{t+1} \mathcal{X}_t,
\end{equation}
where
		$\mathcal{X}_t = \left( \mathcal{A} - \mathcal{B} \mathcal{D}_t^\dagger \mathcal{C}_t \right) - \mathcal{B} \mathcal{F}_t \left( \mathcal{F}^T_t \mathcal{B}^T \mathbf{Q}_{t+1} \mathcal{B} \mathcal{F}_t \right)^\dagger \times \\
		\mathcal{F}_t^T \mathcal{B}^T \mathbf{Q}_{t+1} \left( \mathcal{A} - \mathcal{B} \mathcal{D}_t^\dagger \mathcal{C}_t \right).$
Thus, $\mathcal{U}_t = \mathcal{X}^T_t \mathbf{Q}_{t+1} \mathcal{X}_t \succeq 0$, since $\mathbf{Q}_{t+1} \succeq 0$. 
\end{IEEEproof}

\subsection{Proof of Lemma~\ref{lem: solutionExistence}}\label{sect: solutionExistenceProof}
\begin{IEEEproof}
	Trivially, we have $\mathscr{R} \left(\mathcal{F}_t^T \mathcal{B}^T \mathbf{Q}_{t+1} \mathcal{B} \mathcal{F}_t \right) \subseteq \mathscr{R} \left( \mathcal{F}_t^T \mathcal{B}^T \right).$
We now show that $\mathscr{R} \left( \mathcal{F}_t^T \mathcal{B}^T \right) \subseteq \mathscr{R} \left(\mathcal{F}_t^T \mathcal{B}^T \mathbf{Q}_{t+1} \mathcal{B} \mathcal{F}_t \right)$ by showing that $\mathscr{N} \left(\mathcal{F}_t^T \mathcal{B}^T \mathbf{Q}_{t+1} \mathcal{B} \mathcal{F}_t \right) \subseteq \mathscr{N} \left( \mathcal{B} \mathcal{F}_t \right)$.

Let $\mu \in \mathscr{N} \left(\mathcal{F}_t^T \mathcal{B}^T \mathbf{Q}_{t+1} \mathcal{B} \mathcal{F}_t \right)$, and, by contradiction, suppose that $\mu \notin \mathscr{N} \left( \mathcal{B} \mathcal{F}_t \right)$. Then, we have $\mathcal{F}_t \mu \neq 0$. By definition of $\mathcal{F}_t$, we have $\mathcal{D}_t \mathcal{F}_t \mu = 0,$ which means that $\Psi \mathcal{F}_t {\color{black} \mu} = 0$. Since the matrix $\left[\begin{array} {c} \Gamma \\ \Psi \end{array} \right]$ is injective, $\Psi \mathcal{F}_t \mu = 0$ means that $\widetilde{\mu} = \Gamma \mathcal{F} \mu \neq 0.$
Using the results of Lemma~\ref{lem: QtPSD} and the structure of $\mathcal{H}$ and $\mathcal{B}$, we have
\begin{equation}\label{eqn: solutionExistence1}
	\begin{split}
	\mu^T \mathcal{F}_t^T \mathcal{B}^T \mathbf{Q}_{t+1} \mathcal{B} \mathcal{F}_t \mu =& \mu^T \mathcal{F}_t^T \mathcal{B}^T \mathcal{U}_{t+1} \mathcal{B} \mathcal{F}_t \mu + \\ &\widetilde{\mu}^T Q_{t+1} \widetilde{\mu}.
	\end{split}
\end{equation}
The first term on the right hand side of~\eqref{eqn: solutionExistence1} is nonnegative since $\mathcal{U}_{t+1} \succeq 0$, and the second term on the right hand side of~\eqref{eqn: solutionExistence1} is positive since $\widetilde{\mu} \neq 0$ and $Q_{t+1} \succ 0$. Thus, we have $\mu^T \mathcal{F}_t^T \mathcal{B}^T \mathbf{Q}_{t+1} \mathcal{B} \mathcal{F}_t \mu > 0,$
which contradicts the fact that $\mu \in \mathscr{N} \left(\mathcal{F}_t^T \mathcal{B}^T \mathbf{Q}_{t+1} \mathcal{B} \mathcal{F}_t \right)$. Thus, we have $\mu \in \mathscr{N} \left( \mathcal{B} \mathcal{F}_t \right)$ and $\mathscr{N} \left(\mathcal{F}_t^T \mathcal{B}^T \mathbf{Q}_{t+1} \mathcal{B} \mathcal{F}_t \right) \subseteq \mathscr{N} \left( \mathcal{B} \mathcal{F}_t \right)$. This means that $\mathscr{R} \left( \mathcal{F}_t^T \mathcal{B}^T \right) \subseteq \mathscr{R} \left(\mathcal{F}_t^T \mathcal{B}^T \mathbf{Q}_{t+1} \mathcal{B} \mathcal{F}_t \right)$.
\end{IEEEproof}

\subsection{Proof of Lemma~\ref{lem: zUnique}}\label{sect: zUniqueProof}
\begin{IEEEproof}
	Let $z_1, z_2 \in Z_t(\psi)$, and suppose, by contradiction, that $\mathcal{F}_t z_1 \neq \mathcal{F}_t z_2$. Let $\mu = z_1 - z_2$, which means $\mathcal{F}_t \mu \neq 0$. Since $z_1, z_2 \in Z_t \left( \psi \right)$, we have
$  \mathcal{F}_t^T \mathcal{B}^T \mathbf{Q}_{t+1} \mathcal{B} \mathcal{F}_t \mu = 0. $
Following the proof of Lemma~\ref{lem: solutionExistence}, we have that if $\mathcal{F}_t \mu \neq 0$, then
$\mu^T \mathcal{F}_t^T \mathcal{B}^T \mathbf{Q}_{t+1} \mathcal{B} \mathcal{F}_t \mu > 0,$ which contradicts the fact that $ \mathcal{F}_t^T \mathcal{B}^T \mathbf{Q}_{t+1} \mathcal{B} \mathcal{F}_t \mu = 0.$ Thus, we have $\mathcal{F}_t z_1 = \mathcal{F}_t z_2$.
\end{IEEEproof}

\bibliography{IEEEabrv,References}

\begin{thebibliography}{10}
\providecommand{\url}[1]{#1}
\csname url@samestyle\endcsname
\providecommand{\newblock}{\relax}
\providecommand{\bibinfo}[2]{#2}
\providecommand{\BIBentrySTDinterwordspacing}{\spaceskip=0pt\relax}
\providecommand{\BIBentryALTinterwordstretchfactor}{4}
\providecommand{\BIBentryALTinterwordspacing}{\spaceskip=\fontdimen2\font plus
\BIBentryALTinterwordstretchfactor\fontdimen3\font minus
  \fontdimen4\font\relax}
\providecommand{\BIBforeignlanguage}[2]{{%
\expandafter\ifx\csname l@#1\endcsname\relax
\typeout{** WARNING: IEEEtran.bst: No hyphenation pattern has been}%
\typeout{** loaded for the language `#1'. Using the pattern for}%
\typeout{** the default language instead.}%
\else
\language=\csname l@#1\endcsname
\fi
#2}}
\providecommand{\BIBdecl}{\relax}
\BIBdecl

\bibitem{Cardenas}
A.~A. C\'{a}rdenas, S.~Amin, Z.~Lin, Y.~Huang, C.~Huang, and S.~Sastry,
  ``Attacks against process control systems: Risk assessment, detection, and
  response,'' in \emph{Proc. 6th ACM Symposium on Information, Computer and
  Communications Security}, Hong Kong, Mar. 2011, pp. 355--366.

\bibitem{CardenasOld}
A.~A. C\'{a}rdenas, S.~Amin, and S.~Sastry, ``Research challenges for the
  security of control systems,'' in \emph{Proc. of the 3rd Conf. on Hot Topics
  in Security}, San Jos\'{e}, CA, Jul. 2008, pp. 1--6.

\bibitem{CarAttack}
K.~Koscher, A.~Czeskis, F.~Roesner, S.~Patel, T.~Kohno, S.~Checkoway, D.~McCoy,
  B.~Kantor, D.~Anderson, H.~Shacham, and S.~Savage, ``Experimental security
  analysis of a modern automobile,'' in \emph{Proc. of the 2010 {IEEE}
  Symposium on Security and Privacy}, Oakland, CA, May 2010, pp. 447--462.

\bibitem{Iran}
S.~Peterson and P.~Faramarzi, ``Iran hijacked {US} drone, says {Iranian}
  engineer,'' \emph{Christian Science Monitor}, vol.~15, Dec. 2011.

\bibitem{DroneHack}
D.~Shepard, J.~Bhatti, and T.~Humphreys, ``Drone hack,'' \emph{GPS World},
  vol.~23, no.~8, pp. 30--33, Aug. 2012.

\bibitem{TeixeiraModels}
A.~Teixeira, D.~P\'{e}rez, H.~Sandberg, and K.~H. Johansson, ``Attack models
  and scenarios for networked control systems,'' in \emph{Proc. 1st ACM
  International Conf. on High Confidence Networked Systems}, Beijing, China,
  Apr. 2012, pp. 55--64.

\bibitem{Mo}
Y.~Mo and B.~Sinopoli, ``Integrity attacks on cyber-physical systems,'' in
  \emph{Proc. of the 1st ACM International Conf. on High Confidence Networked
  Systems}, Beijing, China, Apr. 2012, pp. 47--54.

\bibitem{Liu}
Y.~Liu, M.~K. Reiter, and P.~Ning, ``False data injection attacks against power
  systems in electric power grids,'' in \emph{Proc. of the 16th ACM Conf. on
  Computer and Communications Security}, Chicago, IL, Nov. 2009, pp. 21--32.

\bibitem{Fawzi}
H.~Fawzi, P.~Tabuada, and S.~Diggavi, ``Secure estimation and control for
  cyber-physical systems under adversarial attacks,'' \emph{{IEEE} Trans.
  Autom. Control}, vol.~59, no.~6, pp. 1454--1467, Jun. 2014.

\bibitem{Shoukry}
Y.~Shoukry and P.~Tabuada, ``Event-triggered state observers for sparse sensor
  noise/attack,'' \emph{{IEEE} Trans. Autom. Control}, vol.~61, no.~8, pp.
  2079--2091, Aug. 2016.

\bibitem{Pajic}
M.~Pajic, J.~Weimer, N.~Bezzo, P.~Tabuada, O.~Sokolsky, I.~Lee, and G.~Pappas,
  ``Robustness of attack-resilient state estimators,'' in \emph{2014 ACM IEEE
  International Conf. on Cyber-Physical Systems}, Berlin, Germany, Apr. 2004,
  pp. 163--174.

\bibitem{Pasqualetti}
F.~Pasqualetti, F.~Dorfler, and F.~Bullo, ``Attack detection and identification
  in cyber-physical systems,'' \emph{{IEEE} Trans. Autom. Control}, vol.~58,
  no.~11, pp. 2715--2729, Nov. 2013.

\bibitem{ChenICASSP}
Y.~Chen, S.~Kar, and J.~M.~F. Moura, ``{C}yber-{P}hysical {S}ystems: Dynamic
  sensor attacks and strong observability,'' in \emph{Proc. of the 40th IEEE
  International Conf. on Acoustics, Speech and Signal Processing (ICASSP)},
  Brisbane, Australia, Apr. 2015, pp. 1752--1756.

\bibitem{ChenSideInfo}
------, ``Dynamic attack detection in {C}yber-{P}hysical {S}ystems with side
  initial state information,'' \emph{{IEEE} Trans. Autom. Control}, vol.~PP,
  no.~99, pp. 1--6, Nov. 2016.

\bibitem{MoWorkshop}
Y.~Mo and B.~Sinopoli, ``False data injection attacks in control systems,'' in
  \emph{Proc. of the 1st Workshop on Secure Control Systems}, Stockholm,
  Sweden, Apr. 2010, pp. 56--62.

\bibitem{ScalarAttack}
R.~Zhang and P.~Venkitasubramaniam, ``Stealthy control signal attacks in scalar
  lqg systems,'' in \emph{Proc. of the 2015 IEEE Global Conf. on Signal and
  Information Processing (GLOBALSIP)}, Orlando, FL, Dec. 2015, pp. 240--244.

\bibitem{ChenACC}
Y.~Chen, S.~Kar, and J.~M.~F. Moura, ``Cyber-physical attacks constrained by
  control objectives,'' in \emph{Proc. of the 2016 American Control Conf.
  (ACC)}, Boston, MA, Jul. 2016, pp. 1185--1190.

\bibitem{ChenAttack1}
------, ``Cyber physical attacks with control objectives,'' \emph{ArXiv
  e-prints}, Jul. 2016.

\bibitem{Willsky}
A.~S. Willsky, ``A survey of design methods for failure detection in dynamic
  systems,'' \emph{Automatica}, vol.~12, no.~6, pp. 601--611, Nov. 1976.

\bibitem{MoAuthentication}
Y.~Mo, S.~Weerakkody, and B.~Sinopoli, ``Physical authentication of control
  systems: Designing watermarked control inputs to detect counterfeit sensor
  outputs,'' \emph{{IEEE} Control Syst. Mag.}, vol.~35, no.~1, pp. 99--109,
  Feb. 2015.

\bibitem{MoScada}
Y.~Mo, R.~Chabukswar, and B.~Sinopoli, ``Detecting integrity attacks on {SCADA}
  systems,'' \emph{{IEEE} Trans. Control Syst. Technol.}, vol.~22, no.~4, pp.
  1396--1407, Jul. 2014.

\bibitem{RecursiveFeasibility}
J.~L\"{o}fberg, ``Oops! {I} cannot do it again: Testing for recursive
  feasibility in {MPC},'' \emph{Automatica}, vol.~48, no.~3, pp. 550--555, Mar.
  2012.

\bibitem{Rappaport}
D.~Rappaport and L.~M. Silverman, ``Structure and stability of discrete-time
  optimal systems,'' \emph{{IEEE} Trans. Autom. Control}, vol.~16, no.~3, pp.
  227--233, Jun. 1971.

\bibitem{ChenCDC}
Y.~Chen, S.~Kar, and J.~M.~F. Moura, ``Cyber physical attacks with control
  objectives and detection constraints,'' in \emph{Proc. of the 55th IEEE Conf.
  on Decision and Control (CDC)}, Las Vegas, NV, Dec. 2016, pp. 1125--1130.

\bibitem{Speyer}
J.~L. Speyer and W.~H. Chung, \emph{Stochastic Processes, Estimation, and
  Control}.\hskip 1em plus 0.5em minus 0.4em\relax SIAM, 2008, ch.~9.

\bibitem{PajicCDC}
M.~Pajic, P.~Tabuada, I.~Lee, and G.~Pappas, ``Attack-resilient state
  estimation in the presence of noise,'' in \emph{Proc. of the 54th IEEE Conf.
  on Decision and Control (CDC)}, Osaka, Japan, Dec. 2015, pp. 5827--5832.

\bibitem{MoDegradation}
Y.~Mo and B.~Sinopoli, ``On the performance degradation of cyber-physical
  systems under stealthy integrity attacks,'' \emph{{IEEE} Trans. Autom.
  Control}, pp. 1--6, Nov. 2015, {E}arly Access.

\bibitem{Boyd}
S.~Boyd, \emph{Convex Optimization}.\hskip 1em plus 0.5em minus 0.4em\relax
  Cambridge University Press, 2004, ch.~9.

\bibitem{Helicopter}
H.~R. Dharmayanda, A.~Budiyono, and T.~Kang, ``State space identification and
  implementation of ${H}_\infty$ control design for small-scale helicopter,''
  \emph{Aircraft Engineering and Aerospace Technology}, vol.~82, no.~6, pp.
  340--352, 2010.

\bibitem{mosek}
\BIBentryALTinterwordspacing
{MOSEK ApS}, \emph{The MOSEK optimization toolbox for MATLAB manual. Version
  7.1 (Revision 28).}, 2015. [Online]. Available:
  \url{http://docs.mosek.com/7.1/toolbox/index.html}
\BIBentrySTDinterwordspacing

\end{thebibliography}

\begin{IEEEbiography}[{\includegraphics[width=1in,height=1.25in,clip,keepaspectratio]{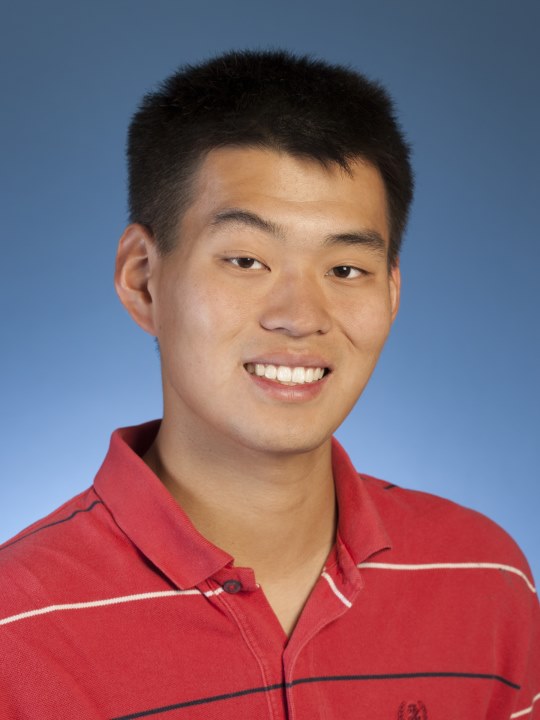}}]{Yuan Chen}(S'14) received the B.S.E. degree in electrical engineering from Princeton University, Princeton, NJ, in 2013. Since 2013, he has been pursuing the Ph.D. degree in the electrical and computer engineering at Carnegie Mellon University, Pittsburgh, PA. His current research activities are focused on distributed inference and the security of cyber-physical systems.
\end{IEEEbiography}

\begin{IEEEbiography}[{\includegraphics[width=1in,height=1.25in,keepaspectratio]{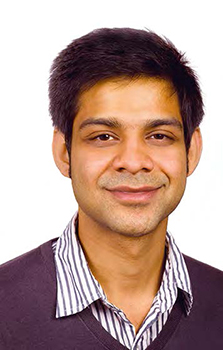}}]{Soummya Kar}
(S'05--M'10)  received a B.Tech. in electronics and electrical communication engineering from the Indian Institute of Technology, Kharagpur, India, in May 2005 and a Ph.D. in electrical and computer engineering from Carnegie Mellon University, Pittsburgh, PA, in 2010. From June 2010 to May 2011, he was with the Electrical Engineering Department, Princeton University, Princeton, NJ, USA, as a Postdoctoral Research Associate.

He is currently an Associate Professor of Electrical and Computer Engineering at Carnegie Mellon University, Pittsburgh, PA, USA. His research interests include decision-making in large-scale networked systems, stochastic systems, multi-agent systems and data science, with applications to cyber-physical systems and smart energy systems. He has published extensively in these topics with more than 140 articles in journals and conference proceedings and holds multiple patents. Recent recognition of his work includes the 2016 O. Hugo Schuck Best Paper Award from the American Automatic Control Council, the 2016 Dean’s Early Career Fellowship from CIT, Carnegie Mellon, and the 2011 A.G. Milnes Award for best PhD thesis in Electrical and Computer Engineering, Carnegie Mellon University.
\end{IEEEbiography}

\begin{IEEEbiography}[{\includegraphics[width=1in,height=1.25in,clip,keepaspectratio]{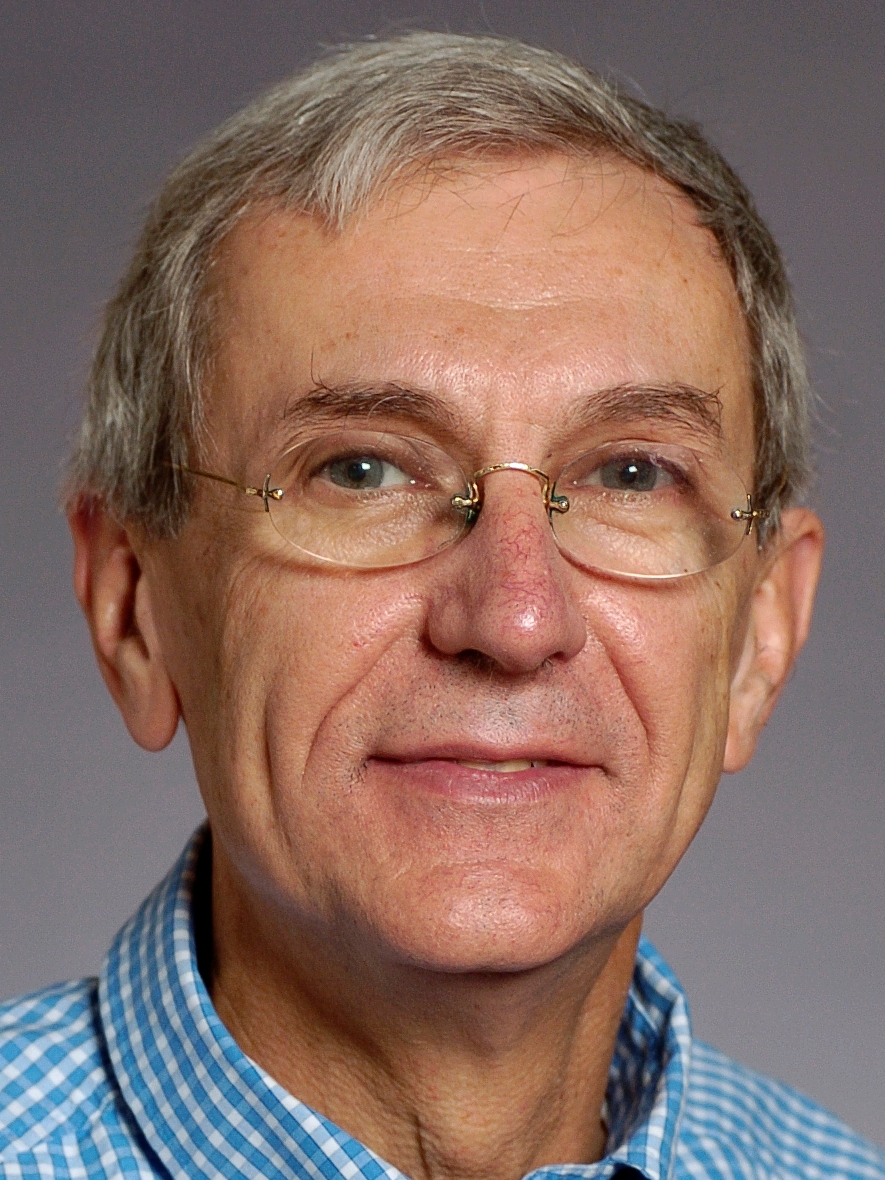}}]{Jos\'e M.~F.~Moura}(S'71--M'75--SM'90--F'94) received the engenheiro electrot\'{e}cnico degree from Instituto Superior T\'ecnico (IST), Lisbon, Portugal, and the M.Sc., E.E., and D.Sc.~degrees in EECS from the Massachusetts Institute of Technology (MIT), Cambridge, MA.

He is the Philip L.~and Marsha Dowd University Professor at Carnegie Mellon University (CMU). He was on the faculty at IST and has held visiting faculty appointments at MIT and New York University (NYU). He founded and directs a large education and research program between CMU and Portugal, www.icti.cmu.edu.

His research interests are on  data science, graph signal processing, and statistical and algebraic signal and image processing. He has published over 550 papers and holds thirteen patents issued by the US Patent Office. The technology of two of his patents (co-inventor A. Kav\v{c}i\'c) are in about three billion disk drives read channel chips of 60~\% of all computers sold in the last 13 years worldwide and was, in 2016, the subject of the largest university verdict/settlement in the information technologies area.

Dr. Moura was the IEEE Technical Activities Vice-President (2016) and member of the IEEE Board of Directors. He served in several other capacities including IEEE Division IX Director, member of several IEEE Boards, President of the IEEE Signal Processing Society(SPS), Editor in Chief for the IEEE Transactions in Signal Processing, interim Editor in Chief for the IEEE Signal Processing Letters.

Dr. Moura has received several awards, including  the Technical Achievement Award and the Society Award from the IEEE Signal Processing Society. In 2016, he received the CMU College of Engineering Distinguished Professor of Engineering Award. He is a Fellow of the IEEE, a Fellow of the American Association for the Advancement of Science (AAAS), a corresponding member of the Academy of Sciences of Portugal, Fellow of the US National Academy of Inventors, and a member of the US National Academy of Engineering.
\end{IEEEbiography} 

\end{document}